\documentclass[sn-mathphys]{sn-jnl}% Math and Physical Sciences Reference Style
\usepackage{graphics}
 \usepackage{graphicx}
\usepackage{amssymb}
\usepackage{amsmath}
\usepackage{color}

\newtheorem{lemma}{Lemma}

\usepackage{hyperref}
\hypersetup{colorlinks=true}
\usepackage{epstopdf}

\jyear{2021}%

\begin{document}

%\title[Product of two continuous random variables]{Product of two continuous random variables – finite and infinite-variance distributions with an application to electricity market transactions}
\title[Product of two continuous random variables]{On the distribution of the product of two continuous random variables with an application to electricity market transactions. Finite and infinite-variance case.}


\author[1]{\fnm{Julia} \sur{Adamska}}\email{249786@student.pwr.edu.pl}

\author[2]{\fnm{{\L}ukasz} \sur{Bielak}}\email{lukasz.bielak@kghm.com}

\author*[1]{\fnm{Joanna} \sur{Janczura}}\email{joanna.janczura@pwr.edu.pl}
\author[1]{\fnm{Agnieszka} \sur{Wy\l oma\'nska}}\email{agnieszka.wylomanska@pwr.edu.pl}

\affil[1]{\orgdiv{Faculty of Pure and Applied Mathematics}, \orgname{Wroclaw University of Science and Technology}, \orgaddress{\street{Wyb. Wyspia\'nskiego 27}, \city{ Wroc\l aw}, \postcode{50-370},  \country{Poland}}}

\affil[2]{\orgname{KGHM}, \orgaddress{\street{M. Sk\l odowskiej-Curie 48}, \city{M. Sk\l odowskiej-Curie 48}, \postcode{59-301}, \country{Poland}}}


\abstract{In this paper we study the distribution of a product of two continuous random variables. We derive formulas for the probability density functions and moments of the products of the Gaussian, log-normal, Student's t and Pareto random variables. In all cases we analyze separately independent as well as correlated random variables. %For the finite-variance case considered in this paper, the zero correlation coefficient indicates the independence of marginal random variables, while for the considered infinite-variance distributions this fact is not true.
%The most interesting case is the product of random variables representing different classes of distributions. 
Based on the theoretical results we use the general maximum likelihood approach for the estimation of the parameters of the product random variables and apply the methodology for a real data case study. We analyze a distribution of the transaction values, being a product of prices and volumes, from a continuous trade on the German intraday electricity market.}

\keywords{product, Gaussian distribution, log-normal distribution, Pareto distribution, Student's t distribution, correlation coefficient, electricity market}


\pacs[MSC Classification]{60E05, 91B26, 62E15, 65C20, 33-02}

\maketitle

\section{Introduction}
One of important problems in statistics and applied mathematics is the recognition of distributions corresponding to various functions of two (or even more) random variables. In the literature different transformations of two random variables, such as a sum, product, ratio, etc., are considered. In this paper we  focus on the product function, i.e we analyze the random variable defined as
\begin{equation}\label{main}
    Z = XY,
\end{equation}
where $X$ and $Y$, called  marginal random variables, have continuous distributions.

The considered problem is not new and there are various examples of random variables for which the product is analyzed in theory and in practice. 
Particular emphasis is placed on the product of two random variables that come from the same class of distributions, see eg. \cite{logistic,pearson,eliptically,trapezoidal}. In this case one can find interesting approaches used in the analysis of the distribution and probabilistic properties of the product random variable. The special attention is paid to the case when two considered random variables are Gaussian or Student's t distributed, see e.g. \cite{Ahsanullah_2014,student_t,gauss2,gauss1,student_statistical}.  The other interesting cases one can find also in \cite{exponential}, where exponentially distributed random variables are considered or in \cite{dirichlet}, where Dirichlet distributed random variables are examined. For other references, see also \cite{beta,beta2,pareto2,beta_statistical}. 
Different classes of distributions are also considered in the literature and the product of such random variables is analyzed, see eg. \cite{normal_laplace,gauss_laplace2,gama_weibul,gama_beta,pareto_gama,beta_statistical} and references therein. 
 %The interesting examples are for instance Gaussian and Laplace distributions \cite{normal_laplace,gauss_laplace2}, Gamma and Weibull distributions \cite{gama_weibul}, Gamma and Beta distributions \cite{gama_beta} or Pareto and Gamma distribution \cite{pareto_gama}. See also \cite{maxwel,pearson_2,pareto4}, for other references.

Beyond theoretical considerations, the random variables that arise as a product (or other functions) of two  variables have found various interesting applications including finance, risk management, economy, but also  physical sciences, reliability theory, hydrology, and many others, see eg. \cite{appl1,appl2,appl3,appl4,appl5,appl6,appl7,appl8,appl10}.

There are known general formulas that can be useful when one considers the random variable defined in Eq. (\ref{main}). For instance, in the general case the probability density function (PDF) of $Z$ is given by \cite{gauss1,rohatgi2015}
\begin{equation}\label{pdf_Z_dependent}
	   f_{Z}(z) = \int_{-\infty}^{\infty} \frac{1}{\mid x\mid} f_{X,Y} \left(x, \frac{z}{x}\right) dx,
\end{equation}
where $f_{X,Y}(x,y)$, $x,y\in \mathbb{R}$ is a two-dimensional PDF of the vector $(X,Y)$.

In this paper we separately study examples of finite-  and infinite-variance distributions of the marginal random variables.  The sample distributions of the first class considered in this paper are the Gaussian and log-normal,  while of the second class - the Student's t \cite{Ahsanullah_2014} and Pareto distributions \cite{Mardia_1962}. We consider two cases. In the first one the marginal random variables come from the same class of distributions, while in the second one - they are representatives of different classes. 

Obviously, the general properties of the product simplify when the marginal random variables are independent. In that case, the corresponding bi-dimensional random vector has the PDF that is a product of the marginal densities, i.e. we have
\begin{equation} \label{pdf_Z_independent}
	   f_{Z}(z) = \int_{-\infty}^{\infty} \frac{1}{\mid x\mid} f_{X} (x) f_{Y} \left( \frac{z}{x}\right) dx.
\end{equation}
where $f_X(x)$ and $f_Y(y)$ are the PDF's of the  marginal distributions corresponding to $X$ and $Y$, respectively. 

The situation is much more complicated when the marginal random variables are statistically dependent. In that case the simplest way to express the dependence between them is the correlation coefficient. However, the correlation coefficient equal to zero does not always correspond to the case when the marginal random variables are statistically independent. From the distributions considered in this paper, this is the case for the Student's t and Pareto ones. Moreover, for the infinite-variance random variables, the correlation coefficient is infinite as well.  Thus, various other approaches might be also considered. One of the ideas utilizes copula as the dependency measure  \cite{appl5,copula2,copula3}. Other dependency measures, adequate for infinite-variance random variables are presented e.g. in \cite{wylomanska2015codifference,slkezak2019codifference,stable,nowicka1}.

The main goal of the paper is to determine probability density functions of the product random variable for the considered distributions and examine, how the probabilistic properties of the marginal random variables influence the distributional properties of the final product. %In the case when the random variables come from the same class of distributions, the special attention is paid to the analysis on how the dependence  between the considered random variables (expressed by means of the correlation coefficient for the finite-variance case) influences the probabilistic properties of their product. 
%It has been shown, that when the marginal random variables are log-normally distributed, the product random variable has also log-normal distribution and this is the only case examined in this paper for which this property is fulfilled. 

The second goal of the paper is the analysis of the distribution of the product random variable when the marginal random variables come from different classes of distributions. The main attention is paid to the analysis on how the parameter responsible for distribution tail influences the PDF of the product random variable.% The other characteristics such as the expected value (if it exists) and the variance (if it exists) are also discussed. The most interesting case seems to be the Gauss-Pareto and the log-normal-Pareto product, for which we obtain an explicit formula for the PDF of the product random variable. It was shown that the resulting distribution in the Gauss-Pareto  case asymptotically behaves like a power-law function.  

The obtained theoretical properties  are  used for the estimation of the product random variable parameters. We use a general idea based on the maximum likelihood algorithm. With the Monte Carlo simulations we confirm the effectiveness of the  procedure. Finally, the methodology is applied to real data from the energy market. We analyze a distribution of the transaction values, being a product of prices and volumes, from a continuous trade on the German intraday electricity market. We show that using the derived product distribution a good fit to the transaction values data can be obtained, being at the same time consistent with the corresponding prices and volumes distributions.

The rest of the paper is organized as follows: First,  we discuss the probabilistic properties of a random variable that arises as a product of random variables coming from two exemplary finite-variance distributions. 
Next, in Section \ref{sec:3} we provide similar considerations for the exemplary random variables from infinite-variance class of distributions. 
In Section \ref{sec:4} we discuss the case when the marginal random variables come from different classes of distributions. The probabilistic properties of the product random variable are calculated and the influence of the parameters of the marginal distributions is discussed. The next part is devoted to the practical aspects of the obtained results. Namely, in Section \ref{sec:5} we propose a general methodology for the estimation of the parameters for the product random variable, while in Section \ref{sec:6} we apply the methodology to a real-data case. We analyze the transactions data from the German electricity market settled on the EPEX energy exchange. Last section concludes the paper and presents other possible applications of the proposed methodology.

%The cumulative distribution function (CDF) of a random variable $Z$ is further in this paper denoted by $F_Z(z)$.

%\textcolor{blue}{
%Under the condition that $X$ and $Y$ have finite expected values and variances, we have the following
%\begin{eqnarray}\label{moments}
%\mathbb{E}(Z)=\int_{-\infty}^{\infty}zf_{Z}(z)dz,~~Var(Z)=\int_{-\infty}^{\infty}z^2f_{Z}(z)dz -\left(\int_{-\infty}^{\infty}zf_{Z}(z)dz\right)^2.
%\end{eqnarray}
%If the random variables $X$ and $Y$ are independent, then the expected value and the variance of $Z$ are expressed in the means of the corresponding characteristics of $X$ and $Y$, namely
%\begin{equation} \label{expected_value}
%    \mathbb{E}(Z) = \mathbb{E}(X) \mathbb{E}(Y),
%\end{equation}
%\begin{equation} \label{variance}
%    Var(Z) = Var(XY) = \mathbb{E}(X^2Y^2) - [\mathbb{E}(X Y)]^2 = \mathbb{E}(X^2)\mathbb{E}(Y^2) - [\mathbb{E}(X) \mathbb{E}(Y)]^2.
%\end{equation}}

\section{Product of  finite-variance distributed random variables}\label{sec:2}
In this section we consider exemplary distributions of $X$ and $Y$ that belong to the finite-variance class, namely the  Gaussian and log-normal ones. 
We recall the known results for the product of two-dimensional Gaussian distribution and derive the corresponding formulas for the log-normal one.
It is worth mentioning, that for these distributions the independence between $X$ and $Y$ is equivalent to the zero valued correlation coefficient and in that case the PDF of the two-dimensional vector $(X,Y)$ is given by Eq. (\ref{pdf_Z_independent}).
\subsection{Gaussian distribution}
The two-dimensional Gaussian distributed random random vector $(X,Y)$ has the following PDF \cite{Roussas_2015} 
\begin{equation}\label{pdf_gaussian_baza}
        f_{X,Y} (x,y) = \frac{ \exp\Bigg\{ -\frac{1}{2(1-\rho^2)} \left[  \frac{(x-\mu_X)^2}{\sigma_X^2} 
        - 2\rho \left(\frac{x-\mu_X}{\sigma_X}\right)\left(\frac{y-\mu_Y}{\sigma_Y}\right)
        +\frac{(y-\mu_Y)^2}{\sigma_Y^2}\right]\Bigg\}}{2 \pi \sigma_X \sigma_Y \sqrt{1-\rho^2}},
\end{equation}
where $x,y\in \mathbb{R}$, $\rho\in (-1,1)$ is the correlation coefficient between random variables $X$ and $Y$; $\mu_X,\mu_Y\in \mathbb{R}$ are the corresponding expected values, while $\sigma_X, \sigma_Y>0$ are the corresponding standard deviations. The cases $\rho=1$ or $\rho=-1$ are not considered in this paper. However, they correspond to the situation when $Y=aX$, with $a\in \mathbb{R}$. As a consequence, the product random variable $Z$ has a chi-square distribution with one degree of freedom. It is easy to see that, when $\rho=0$, the PDF of the random vector $(X,Y)$ is just a product of the PDFs of the Gaussian distributed random variables.

For simplicity we assume that $\mu_X=\mu_Y=0$. In that case the PDF of the random variable $Z$ defined in (\ref{main}) is given by \cite{Gaunt_2018}
\begin{equation} \label{gauss2}
        f_Z (z) = \frac{1}{ \pi\sigma_X \sigma_Y \sqrt{1-\rho^2}} \exp\Bigg\{ \frac{\rho z}{\sigma_X \sigma_Y (1-\rho^2)}\Bigg\} K_0\left(\frac{\mid z\mid}{\sigma_X \sigma_Y (1-\rho^2)}\right),
\end{equation}
where  $z\in \mathbb{R}$, $K_0(\cdot)$ is the Bessel function of the second kind with a purely imaginary argument of zero order. Recall that the family of the modified Bessel functions of the second kind for order $v$ is given by $K_v(x) = \int_0^{\infty} \exp\{-x \cosh(t)\}\cosh(vt) dt$. {The formula for $f_Z(\cdot)$ in the general case with $\mu_X=\mu_Y\in \mathbb{R}$ can be found in \cite{general_normal}.}

%One can identify the PDF given in Eq. (\ref{gauss2}) as the PDF of the variance-gamma distribution. Recall that, if a random variable $V$ has the variance-gamma distribution $V \sim VG(r,\theta ,\widetilde{\sigma},\widetilde{\mu})$, where $r>0, \theta \in \mathbb{R}, \widetilde{\sigma}>0, \widetilde{\mu}\in\mathbb{R}$, then its PDF has the following form \cite{Gaunt_2018}
%\begin{equation}
%    f_V(v) = \frac{\exp\Big\{ \frac{\theta}{\widetilde{\sigma}^2}(v-\widetilde{\mu})\Big\}}{\widetilde{\sigma} \sqrt{\pi} \Gamma{(\frac{r}{2})} }
%    \left( \frac{\mid v-\widetilde{\mu}\mid}{2\sqrt{\theta^2 + \widetilde{\sigma}^2}}\right)^{\frac{r-1}{2}}
%    K_{\frac{r-1}{2}} \left(\frac{\sqrt{\theta^2+\widetilde{\sigma}^2}}{\widetilde{\sigma}^2}\mid v-\widetilde{\mu}\mid\right),
%\end{equation}
%where $v\in \mathbb{R}$, $K_v(\cdot)$ is given by (\ref{bessel}). 
One can show that the PDF given in Eq. (\ref{gauss2}) corresponds to the  variance-gamma distribution PDF: $VG(1, \rho\sigma_x\sigma_y, \sigma_x\sigma_y\sqrt{1-\rho^2},0)$.
%When $\rho=0$ the PDF of $Z$  takes the form \cite{Ahsanullah_2014}
%\begin{equation}
%        f_Z (z) = \frac{1 }{ \sigma_X\sigma_Y\pi} K_0 \left( \frac{\mid z\mid}{\sigma_X\sigma_Y} \right).
%\end{equation}
The expected value and the variance of $Z$ in the general case are given by \cite{Craig_1936}
\begin{align}
    \mathbb{E}(Z) &= \mu_X\mu_Y + \rho \sigma_X\sigma_Y,\\
    Var(Z) &= \sigma_X^2 \mu_Y^2 + \sigma_Y^2\mu_X^2+\sigma_X^2\sigma_Y^2  + 2\rho\mu_X\mu_Y\sigma_X\sigma_Y+\rho^2\sigma_X^2\sigma_Y^2.
\end{align}
%what reduces to the following formulas when $\rho=0$ \cite{Craig_1942}
%\begin{equation}
%    \mathbb{E}(Z) = \mu_X\mu_Y,~~   Var(Z) = \sigma_X^2 \mu_Y^2 + \sigma_Y^2\mu_X^2+\sigma_X^2\sigma_Y^2.
%\end{equation}
It is interesting to note that, for correlated marginal variables the individual scale parameters, $\sigma_X, \sigma_Y$ influence also the expected value of the product, $\mathbb{E}(Z)$, and for central marginal distributions we have  $\mathbb{E}(Z)\neq0$. On the other hand, if $\rho=0$, then the product expectation is just the product of the individual means, $\mu_X, \mu_Y$. The variance of the product is equal to the product of individual variances only if $\mu_X=\mu_Y=0$ and the marginal variables are uncorrelated. Otherwise, it is a function of all individual parameters.

%In Fig.  \ref{fig:normal} we demonstrate the PDF  and the distribution tails of the random variable $Z$ that is a product of two central Gaussian distributed random variables with $\mu_X=\mu_Y=0$, $\sigma_X=\sigma_Y=1$ for different values of the $\rho$ parameter.
%The resulting PDF has a characteristic spike at $z=0$. It is symmetric only if $\rho=0$, i.e. the marginal variables are independent. Otherwise the product distribution is right-skewed for $\rho>0$ and left-skewed for $\rho<0$. 

%\begin{figure}[ht]
%    \centering
%    \includegraphics[scale = 0.5]{normal_fz.eps}
%    \caption{Distribution of the product of two Gaussian distributed random variables with $\mu_X = \mu_Y = 0$, $\sigma_X=\sigma_Y=1$ and different values of the $\rho$ parameter. Left panel: the PDF of the random variable $Z$. Right panel: the distribution tail for the random variable $Z$ (in log-log scale). }
%    \label{fig:normal}
%\end{figure}

\subsubsection{Log-normal distribution}
 The two-dimensional log-normally distributed random vector $(X,Y)$ with parameters $\mu_X,\mu_Y\in \mathbb{R}$, $\sigma_X,\sigma_Y>0$ and $\rho\in (-1,1)$ is defined in the following way \cite{log_multi}
\begin{eqnarray}\label{rep1}
(X,Y)=\exp\{(N_1,N_2)\},
\end{eqnarray}
where $(N_1,N_2)$ is the Gaussian random vector defined by the PDF in Eq. (\ref{pdf_gaussian_baza}) with the parameters $\mu_X,\mu_Y,\sigma_X,\sigma_Y$ and $\rho$. Thus, the marginal random variables $X$ and $Y$ have the representations
\begin{eqnarray}\label{rep2}
X=\exp\{N_1\}, ~Y=\exp\{N_2\}
\end{eqnarray}
and in the general case the $\rho$ parameter is the correlation between $N_1$ and $N_2$. This influences also the correlation between $X$ and $Y$. Indeed, one can easily show that the covariance between $X$ and $Y$ is as follows
\begin{eqnarray}
\text{cov}(X,Y)=\exp\left\{\mu_X+\mu_Y+\frac{1}{2}(\sigma_X^2+\sigma_Y^2)\right\}\Bigg[\exp\left\{\rho\sigma_X\sigma_Y\right\}-1\Bigg].
\end{eqnarray}
%\begin{eqnarray}
%\text{corr}(X,Y)=\frac{\text{cov}(X,Y)}{\exp\Bigg\{2\mu_X + \sigma_X^2+2\mu_Y + \sigma_Y^2\Bigg\}\left[  \exp\Bigg\{ \sigma_X^2 + \sigma_Y^2\Bigg\} -1 \right]}.
%\end{eqnarray}
The PDF of $(X,Y)$ is given by  \cite{Yerel_BVLN}
\begin{equation}\label{log-normal}
        f_{X,Y} (x,y) = \frac{\exp\Bigg\{ -\frac{\frac{(\log(x)-\mu_X)^2}{\sigma_X^2} 
        - 2\rho \frac{(\log(x)-\mu_X)(\log(y)-\mu_Y)}{\sigma_X \sigma_Y}
        +\frac{(\log(y)-\mu_Y)^2}{\sigma_Y^2}}{2(1-\rho^2)} \Bigg\}}{2 \pi x y \sigma_X \sigma_Y \sqrt{1-\rho^2}} 
        ,
\end{equation}
where $x,y >0$. If $\rho=0$, the marginal variables $X$ and $Y$ are independent and the PDF given in (\ref{log-normal}) is a product of PDFs of one-dimensional random variables with the  log-normal distribution. The cases $\rho=1$ and $\rho=-1$ are not considered in this paper. However, one can notice that in these cases  the product random variable $Z$ has still the  log-normal distribution.

One can show that the random variable $Z$ defined in Eq. (\ref{main}) in the considered case has also one-dimensional log-normal distribution. Indeed, from Eq. (\ref{rep2}) one has
\[Z=XY=\exp\{N_1+N_2\}.\]
Since $N_1$ and $N_2$ are jointly Gaussian, the random variable $N_1+N_2$ still has the Gaussian distribution with the expected value $\mu_X+\mu_Y$ and the variance $\sigma_X^2+\sigma_Y^2+2\rho\sigma_X\sigma_Y$. Thus, the PDF of $Z$ has the following form\begin{eqnarray}
f_Z(z)=\frac{1}{\sqrt{2\pi(\sigma_X^2+\sigma_Y^2+2\rho\sigma_X\sigma_Y)}z}\exp\left\{-\frac{1}{2}\frac{(\log(z)-\mu_X-\mu_Y)^2}{\sigma_X^2+\sigma_Y^2+2\rho\sigma_X\sigma_Y}\right\}
\end{eqnarray}
for $z>0$. %If $\rho=0$, then we have
%\begin{eqnarray}
%f_Z(z)=\frac{1}{\sqrt{2\pi(\sigma_X^2+\sigma_Y^2)}z}\exp\left\{-\frac{1}{2}\frac{(\log(z)-\mu_X-\mu_Y)^2}{\sigma_X^2+\sigma_Y^2}\right\}, ~~z>0.
%\end{eqnarray}
%\textcolor{blue}{One can easily show that
%\begin{eqnarray}\label{Gaussian_moments}
%\mathbb{E}(Z)&=&\exp\left\{\mu_X+\mu_Y+\frac{1}{2}(\sigma_X^2+\sigma_Y^2+2\rho\sigma_X\sigma_Y)\right\}, \\
%Var(Z)&=&\exp\Bigg\{2(\mu_X+\mu_Y)+\sigma_X^2+\sigma_Y^2+2\rho\sigma_X\sigma_Y\Bigg\}\nonumber \\
%&& \times \Bigg[\exp\left\{\sigma_X^2+\sigma_Y^2+2\rho\sigma_X\sigma_Y\right\}-1\Bigg].
%\end{eqnarray}}
%what in the case of $\rho=0$ reduces to the following formulas
%\begin{align}
%    \mathbb{E}(Z) &= \exp\Bigg\{\mu_X + \mu_Y + \frac{\sigma_X^2 + %\sigma_Y^2}{2}\Bigg\},\\
 %Var(Z) &= \exp\Bigg\{2\mu_X + \sigma_X^2+2\mu_Y + \sigma_Y^2\Bigg\}\left[  %\exp\Bigg\{ \sigma_X^2 + \sigma_Y^2\Bigg\} -1 \right].
%\end{align}
The resulting product distribution is still log-normal with the following parameters: $\mu_Z=\mu_X+\mu_Y$ and $\sigma^2_Z=\sigma_X^2+\sigma_Y^2+2\rho\sigma_X\sigma_Y$. Hence, the $\rho$ coefficient influences the scale parameter of the resulting distribution. %It is illustrated in Fig.  \ref{fig:lognormal} were we plot the PDF of the random variable $Z$ and the corresponding distribution tail for different values of the $\rho$ parameter. As can be observed, higher values of $\rho$ result in lower probabilities of extreme observations.
%\begin{figure}[ht]
%    \centering
%    \includegraphics[scale = 0.5]{lognormal_fz.eps}
%    \caption{Distribution of the product of two log-normally distributed random variables with $\mu_X = \mu_Y = 0$, $\sigma_X=\sigma_Y=1$ and different values of the $\rho$ parameter. Left panel: the PDF of the  random variable $Z$. Right panel: the distribution tail of the random variable $Z$ (in log-log scale).}
%    \label{fig:lognormal}
%\end{figure}

\section{Product of infinite-variance distributed random variables}  \label{sec:3}
In this section, we consider exemplary distributions of $X$ and $Y$ that belong to the infinite-variance class of distributions. The considered cases are the Student's t and Pareto. In contrast to the finite-variance  distributions considered in the previous section, the zero correlation coefficient does not correspond to the case when the marginal random variables are statistically independent.
\subsection{Student's t distribution}
A common method of construction of a two-dimensional Student's t distributed random vector $(X,Y)$ is based on the following observation. Let us assume  that $(N_1,N_2)$ is the two-dimensional Gaussian vector defined by the PDF  in Eq. (\ref{pdf_gaussian_baza})  with the parameters $\mu_X=\mu_Y=0$, $\sigma_X=\sigma_Y=1$ and $\rho\in (-1,1)$, and $\chi^2$ is the one-dimensional random variable with chi-square distribution with $n>0$ degrees of freedom. Moreover, assume that $(N_1,N_2)$ and $\chi^2$ are independent. Then the random vector defined as
\begin{eqnarray}
(X,Y)=\frac{(N_1,N_2)}{\sqrt{\chi^2/n}}
\end{eqnarray}
has a two-dimensional Student's t distribution with $n$ degrees of freedom and its PDF is given by \cite{Lai_2009}
\begin{equation}\label{student1}
    f_{X,Y}(x,y) = \frac{1}{2\pi \sqrt{1-\rho^2}}\left[ 1 + \frac{x^2 - 2\rho xy + y^2}{n(1-\rho^2)} \right]^{-\frac{n+2}{2}},~~x,y\in\mathbb{R}.
\end{equation}
The marginal random variable  $X$ (and $Y$) has the one-dimensional Student's t distribution defined by the PDF given by \cite{student_baza}
\begin{eqnarray}\label{student2}
f_X(x)=\frac{\Gamma((n+1)/2)}{\sqrt{n\pi}\Gamma(n/2)}\left(1+\frac{x^2}{n}\right)^{-(n+1)/2},~x\in\mathbb{R},
\end{eqnarray}
where $\Gamma(\cdot)$ is the gamma function, i.e. $\Gamma(\alpha) = \int_{0}^{\infty} t^{\alpha-1} e^{-t} dt$ for $\alpha$ such that $\text{Re}(\alpha)>0$. Note that the number of degrees of freedom, $n$, is equal for both marginal variables. %The expected value (if $n>1$) and the variance (if $n>2$) of the marginal random variable $X$ are given by
%\begin{eqnarray}
%\mathbb{E}(X)=0,~Var(X)=\frac{n}{n-2}.
%\end{eqnarray}
For $n>2$ the random variables $X$ and $Y$ described by the PDF in (\ref{student1}) have the following covariance  
\begin{eqnarray}
\text{cov}(X,Y)=\frac{n}{n-2}\rho.%,~\text{corr}(X,Y)=\rho.
\end{eqnarray}
One can see that, although the zero correlation coefficient  corresponds to $\text{cov}(X,Y)=0$, it is not equivalent to the independence of the random variables $X$ and $Y$, since in that case the PDF of a random vector $(X,Y)$ (see Eq. (\ref{student1})) is not a product of the PDFs of the marginal distributions (see Eq. \ref{student2}). Similarly, as for the previously considered examples, the cases $\rho=1$ and $\rho=-1$ are not analyzed in this paper.

In the case when the random vector $(X,Y)$ is described by the PDF given in Eq. (\ref{student1}), then using formula (\ref{pdf_Z_dependent}) one obtains 
\begin{eqnarray}\label{student3}
f_{Z}(z) = \frac{1}{2 \pi \sqrt{1-\rho^2}} \int_{-\infty}^{\infty} \frac{1}{\mid x\mid} 
\left[  1 + \frac{x^2 - 2\rho z + \frac{z^2}{x^2}}{n(1-\rho^2)} \right]^{-\frac{n+2}{2}} dx, ~x\in\mathbb{R}.
\end{eqnarray}
The integral in Eq. (\ref{student3}) has no closed form and can be  expressed by means of the Appel hypergeometric special function. However, one can obtain its value with numerical calculations.

If the random variables $X$ and $Y$ are independent and are described by the Student's t distributions of PDFs given in Eq. (\ref{student2}) with degrees of freedom $n_X>0$ and $n_Y>0$ for X and Y distributions, respectively, then  the PDF of the vector $(X,Y)$ is given by a product of marginal densities of $X$ and $Y$. 
%\begin{equation}\label{student4}
%f_{X,Y}(x,y)=\frac{\Gamma\left(\frac{n_X+1}{2}\right)\Gamma\left(\frac{n_y+1}{2}\right)}{\sqrt{n_Xn_y}\pi\Gamma\left(\frac{n_y}{2}\right)\Gamma\left(\frac{n_X}{2}\right)}\left(1+\frac{x^2}{n_X}\right)^{-(n_X+1)/2}\left(1+\frac{y^2}{n_y}\right)^{-(n_y+1)/2}.
%\end{equation}
Thus, the product random variable $Z$ has the following PDF for $z\in \mathbb{R}$ \cite{Ahsanullah_2014}
\begin{equation} \label{student4}
f_{Z}(z) =\frac{2\Gamma{\left(\frac{n_X+1}{2}\right)}\Gamma{\left(\frac{n_Y+1}{2}\right)}}{\sqrt{ n_X n_Y } \pi\Gamma{\left(\frac{n_X}{2}\right)}\Gamma{\left(\frac{n_Y}{2}\right)}}   \int_{0}^{\infty} \frac{1}{x} 
\left( 1+\frac{x^2}{n_X}  \right)^{-\frac{n_X+1}{2}}
\left( 1+\frac{z^2}{x^2 n_Y}  \right)^{-\frac{n_Y+1}{2}}
dx.
\end{equation}
Similarly as in the previous case, the above integral can only be calculated numerically, as it has no closed form representation. 

When one considers the random vector of the two-dimensional Student's t distribution defined by the PDF in Eq. (\ref{student1}), then the expectation (when $n>1$) and the variance (when $n>2$) of the random variable $Z$ can be calculated based on its PDF (see Eq. (\ref{student3})). %formulas given in Eq. (\ref{moments}). 
In this case, there are no closed forms for the mentioned statistics. However, when $X$ and $Y$ are independent random variables, %i.e. when the PDF of a random vector $(X,Y)$ is given by Eq. (\ref{student4}), 
then the expected value (if $n_X,n_y>1$) and the variance (if $n_X,n_Y>2$) of the random variable $Z$ are given by \cite{Ahsanullah_2014}
\begin{eqnarray}
 \mathbb{E}(Z) = 0,~~Var(Z) = \frac{n_X n_Y}{(n_X-2)(n_Y-2)}.
\end{eqnarray}
The variance of the product of independent marginal variables  is just a product of the individual variances. Moreover, it tends to the variance of one of the variables if the number of degrees of freedom  of the second variable, $n_X$ or $n_Y$, goes to infinity. If both parameters, $n_X,n_Y\rightarrow \infty$, then the variance of the product decreases to 1. 

\begin{figure}[ht]
    \centering
    \includegraphics[scale = 0.5]{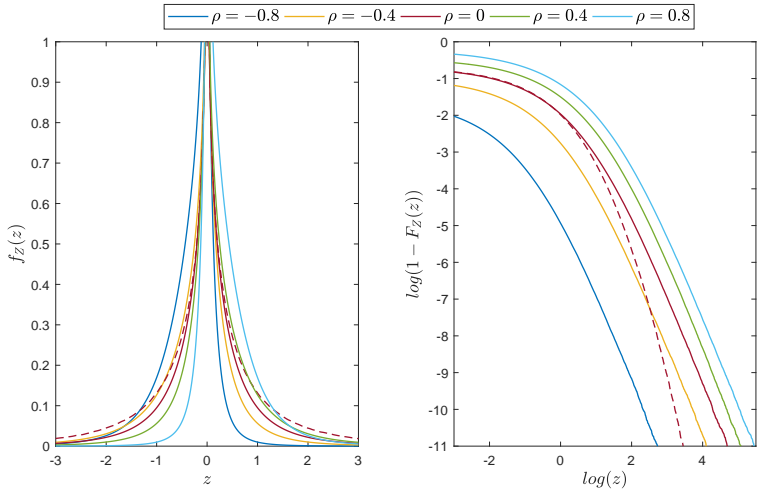}
    \caption{ Distribution of the product of the marginal random variables from the two-dimensional Student's t distribution with $n=5$ and different values of the $\rho$ parameter. Left panel: the PDF of the random variable $Z$. Right panel: the distribution tail for the random variable $Z$ (in log-log scale). The dotted line corresponds to the case when $X$ and $Y$ random variables are independent and $n_X=n_Y=5$. The colors of the solid and dotted lines correspond to the same $\rho$ parameter.}
    \label{fig:student}
\end{figure}

In Fig.  \ref{fig:student} we demonstrate the PDF and the distribution tails (1-cumulative distribution functions, 1-CDFs) of the random variable $Z$ that is a product of marginal random variables from the two-dimensional Student's t distribution described by the PDF (\ref{student1}) for different values of the $\rho$ parameter. %Similarly, as for the Gaussian case (see Fig. \ref{fig:normal}) the 
The resulting PDF is symmetric only if $\rho=0$. Otherwise, it is right-skewed for $\rho>0$ and left-skewed for $\rho<0$. The distribution tails are clearly heavier than in the Gaussian case. For the comparison, we also show the plots for the corresponding product of independent Student's t marginal variables, i.e., with PDF given by Eq. (\ref{student4}). The resulting distribution is again symmetric, however, its tail is lighter than in the corresponding case of dependent marginal variables.

\subsection{Pareto distribution}
The two-dimensional Pareto distributed random vector $(X,Y)$ with parameters $a, \theta_X,\theta_Y>0$ has the following PDF \cite{Mardia_1962}
\begin{equation}\label{pareto1}
    f_{X,Y}(x,y) = \frac{a (a+1) (\theta_X\theta_Y)^{a+1}}{(x\theta_Y + y \theta_X - \theta_X\theta_Y)^{a+2}}, ~x>\theta_X,~y>\theta_Y.
\end{equation}
The above parametrization corresponds to the Pareto distribution of the first kind \cite{paret}. The $a$ parameter is responsible for the behavior of the distribution tail. The marginal distributions of $X$ and $Y$ are given by 
\begin{eqnarray}\label{pareto2}
f_X(x)=\frac{a\theta_X^{a}}{x^{a+1}},~~x>\theta_X, ~~~~f_Y(y)=\frac{a\theta_Y^{a}}{y^{a+1}},~~y>\theta_Y.
\end{eqnarray}
%The expected value (if $a>1$) and variance (if $a>2$) of the marginal random variable $X$ are given by
%\begin{eqnarray}
%\mathbb{E}(X)=\frac{a\theta_X}{a-1},~~Var(X)=\frac{a\theta_X^2}{(a-1)^2(a-2)}.
%\end{eqnarray}
%The analogous representation holds for the random variable $Y$ ( $\theta_X$ is replaced by $\theta_Y$). 
Note that the shape parameter $a$ is the same for both marginal distributions. For $a> 2$, the random variables $X$ and $Y$ described by the PDF (\ref{pareto1}) are positively correlated. The covariance and correlation are given by, respectively
\begin{eqnarray}
\text{cov}(X,Y)=\frac{\theta_X\theta_Y}{(a-1)^2(a-2)^2},~\text{corr}(X,Y)=\frac{1}{a}.
\end{eqnarray}
It is interesting to note that here the correlation between the marginal variables is governed by the shape parameter $a$. As a consequence, with such parametrization, $X$ and $Y$ are always positively correlated.  

When the random vector $(X,Y)$ is described by the PDF given in Eq. (\ref{pareto1}) using  formula (\ref{pdf_Z_dependent}) one obtains
\begin{equation}\label{pareto3}
    f_Z(z) =(a+1)a (\theta_X\theta_Y)^{a+1}  \int\limits_{\theta_X}^{\frac{z}{\theta_Y}} \frac{1}{x}
    (\theta_Y x + \theta_X \frac{z}{x} - \theta_X\theta_Y)^{-a-2} dx,~~z>\theta_X\theta_Y.
 \end{equation}
The integral in Eq. (\ref{pareto3}) has no closed form and can be  expressed by means of the Appel hypergeometric special function. However, one can obtain its value with numerical calculations.

As one can see, the two-dimensional Pareto distribution defined by the PDF in Eq. (\ref{pareto1}) does not cover the case with uncorrelated marginal random variables.  In the case when $X$ and $Y$ are independent and are described by the Pareto distributions of the PDFs given in (\ref{pareto2}) with $a_X,\theta_X>0$ and $a_Y,\theta_Y>0$, respectively, then the PDF of the vector $(X,Y)$ is given by the product of the marginal PDFs, $f_X(x)$ and $f_Y(y)$. 
%\begin{eqnarray}\label{pareto4}
%f_{X,Y}(x,y)=\frac{a_Xa_Y\theta_X^{a_X}\theta_Y^{a_Y}}{x^{a_X+1}y^{a_Y+1}},~x>\theta_X,y>\theta_Y.
%\end{eqnarray}
\begin{lemma}
Assume that the independent random variables $X$ and $Y$ have Pareto distribution with parameters $a_X,\theta_X$ and $a_Y,\theta_Y$, respectively. Then, for $a_X\neq a_Y$, the random variable $Z$ defined in (\ref{main}) has the following PDF 
\begin{equation}\label{pareto_independent}
    f_{Z}(z) = a_X a_Y \left( \frac{z^{-a_X-1} (\theta_X\theta_Y)^{a_X}}{a_Y -a_X}  +  \frac{ z^{-a_Y-1} (\theta_X\theta_Y)^{a_Y}}{a_X - a_Y}  \right),~~  z>\theta_X\theta_Y.
\end{equation}
When $a_X=a_Y=a$ the PDF of $Z$ has the following form
\begin{eqnarray}
f_Z(z)=a^2(\theta_X\theta_Y)z^{-a-1}\log\left(\frac{z}{\theta_X\theta_Y}\right),~~z>\theta_X\theta_Y.
\end{eqnarray}
\end{lemma}
\textbf{Proof:} The proof comes directly from formula (\ref{pdf_Z_independent}) and the fact that the PDF of the random vector $(X,Y)$ is given by the product of marginal PDFs of $X$ and $Y$ 
%in Eq. (\ref{pareto4}) 
for $a_X\neq a_Y$ and $a_X=a_Y=a$. $\Box$.

When one considers a random vector of the two-dimensional Pareto distribution defined by the PDF in Eq. (\ref{pareto1}), then there are no closed forms for the expectation (when $a>1$) and the variance (when $a>2$) of the random variable $Z$% can be calculated using formulas given in Eq. (\ref{moments}). In this case, there are no closed forms for the mentioned statistics
. However, when $X$ and $Y$ are independent random variables, %i.e., when the PDF of a random vector is given by Eq. (\ref{pareto4}), 
then the expected value and the variance of the random variable $Z$ are given by \cite{Mardia_1962} 
\begin{align}
    \mathbb{E}(Z) &= \frac{\theta_X \theta_Y a_X a_Y}{(a_X-1)(a_Y-1)},~~a_X,a_Y>1,\\
    Var(Z)&=  \theta_X^2 \theta_Y^2 a_Xa_Y
    \left[\frac{(a_X-1)^2+(a_Y-1)^2-1}{(a_X-2)(a_Y-2)(a_X-1)^2(a_Y-1)^2}
    \right], ~~a_X,a_Y>2.
\end{align}

\begin{figure}[ht]
    \centering
    \includegraphics[scale = 0.5]{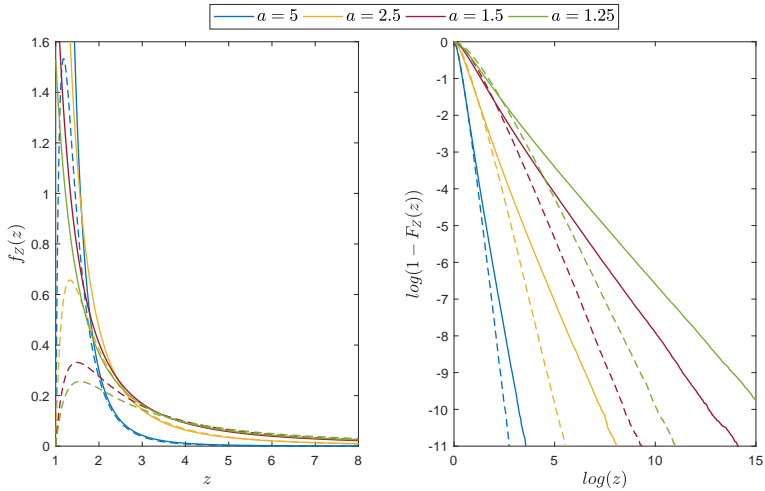}
    \caption{Distribution of the product of marginal random variables from the two-dimensional Pareto distribution with $\theta_X=\theta_Y=1$ and different values of the $a$ parameter. Left panel: the PDF of the random variable $Z$. Right panel: the distribution tail for the random variable $Z$ (in log-log scale). The dotted lines correspond to the case when $X$ and $Y$ are independent and $a_X=a_Y=a$, $\theta_X=\theta_Y=1$. The colors of the solid and dotted lines correspond to the same $a$ parameter.}
    \label{fig:pareto}
\end{figure}

In Fig.  \ref{fig:pareto} we plot the PDF and the corresponding distribution tail of the random variable $Z$ that is a product of  marginal random variables from the two-dimensional Pareto distribution described by the PDF (\ref{pareto1}) for different values of the $a$ parameter. Recall that the correlation (if exists) is directly related to the shape parameter $a$, i.e. $\text{corr}(X,Y)=\frac{1}{a}$. So it is influenced by the distribution tails. With larger $a$ (lower $\text{corr}(X,Y)$) the probability of extreme observations becomes lower. In all cases, the distribution tails behave like a power function with a clear linear shape in the double logarithmic scale. 
For comparison, we also show the plots for the corresponding product of independent Pareto marginal variables%, i.e. with the PDF given by Eq. (\ref{pareto4})
. In this case, the probabilities of extreme observations are lower than for the corresponding dependent Pareto distribution, but the power tail behavior is preserved. 

\section{Product of finite- and infinite-variance distributed random variables} \label{sec:4}
The distribution of the  random variable that arises as a product of two other random variables belonging to different families of distributions is of considerable
importance and current interest in various applications.  In this section, we examine the case when one of the marginal random variables is finite-variance distributed, while the second one belongs to the infinite-variance class. The product random variable shares some properties of both marginal variables, but on the other hand it exhibits additional interesting features. In this section, we assume that the marginal random variables are independent. 
\subsection{Gaussian and Student’s t distributions}
We assume that the one-dimensional random variable $X$ has the Gaussian distribution with parameters $\mu_X\in \mathbb{R}$ and $\sigma_X>0$, while the one-dimensional random variable $Y$ has the Student's t distribution with the PDF given in Eq. (\ref{student2}) with $n_Y>0$.  We assume that $X$ and $Y$ are independent, thus the PDF of the random vector $(X,Y)$ is just the product of the corresponding marginal PDFs of $X$ and $Y$.
%has the following PDF
%\begin{equation} \label{Gauss_Student}
%   f_{X,Y}(x,y) = \frac{\Gamma{\left(\frac{n_Y+1}{2}\right)}}{\pi \sigma_X\sqrt{2n_Y} \Gamma{\left(\frac{n_Y}{2}\right)}} \exp\Bigg\{-\frac{(x-\mu_X)^2}{2\sigma_X^2}\Bigg\}\left( 1+\frac{y^2}{n_Y}  \right)^{-\frac{n_Y+1}{2}}, ~~x,y\in \mathbb{R}.
%\end{equation}
Using Eq. (\ref{pdf_Z_independent}) one obtains the PDF of the random variable $Z$
\begin{eqnarray}\label{GS2}
f_Z(z)&=&\int_{-\infty}^{\infty}\frac{1}{\mid x\mid} f_{X,Y}\left(x,\frac{z}{x}\right)dx\\
&=&\frac{\Gamma{\left(\frac{n_Y+1}{2}\right)}}{\pi \sigma_X\sqrt{2n_Y} \Gamma{\left(\frac{n_Y}{2}\right)}}\int_{-\infty}^{\infty}\frac{1}{\mid x\mid}\exp\Bigg\{-\frac{(x-\mu_X)^2}{2\sigma_X^2}\Bigg\}\left( 1+\frac{z^2}{n_Yx^2}  \right)^{-\frac{n_Y+1}{2}}dx\nonumber\\
&=&\frac{2\Gamma{\left(\frac{n_Y+1}{2}\right)}}{\pi \sigma_X\sqrt{2n_Y} \Gamma{\left(\frac{n_Y}{2}\right)}}\int_{0}^{\infty}\frac{1}{x}\exp\Bigg\{-\frac{(x-\mu_X)^2}{2\sigma_X^2}\Bigg\}\left( 1+\frac{z^2}{n_Yx^2}  \right)^{-\frac{n_Y+1}{2}}dx.\nonumber
\end{eqnarray}
The PDF of $Z$ has no closed form representation and its calculation requires numerical computation of the integral in Eq. (\ref{GS2}).
%In that case the PDF of $Z$ defined in (\ref{main}) is given by, \cite{Ahsanullah_2014}
%\begin{equation} 
 %  f_Z(z) = \frac{\Gamma(\frac{n_Y}{2},\frac{1}{2})}{\sqrt{2n_Y}B(\frac{n_Y}{2},\frac{1}{2})} \psi\left(\frac{n_Y}{2}+\frac{1}{2},1;\frac{z^2}{2\sigma_x^2 n}\right),~z\in\mathbb{R},
%\end{equation}
%where  $\psi(\cdot,\cdot;\cdot)$ is a Kummer's hypergeometric function defined as \textcolor{bulue}{cyt}
 %\begin{equation}
  %   \psi(\alpha,\gamma;z) = \frac{1}{\Gamma(\alpha)} \int_0^{\infty} \exp\{-zt\} t^{\alpha-1}(1+t)^{\gamma-\alpha-1} dt,
 %\end{equation}
%where $Re(\alpha)>0$, $Re(z)>0$. Moreover, for $\alpha>0$ and  $p,q>0$ the functions $\Gamma(\cdot,\cdot)$ and $B(\cdot,\cdot) $ are lower incomplete Gamma and Beta functions, respectively and are defined as, \cite{Ahsanullah_2014}
%\begin{equation}\label{inc}
 %  \Gamma(\alpha,x) = \int_{x}^{\infty} t^{\alpha-1} \exp\{-t\} dt,~~B(p,q) = \frac{\Gamma(p)\Gamma(q)}{\Gamma(p+q)}.
%\end{equation}
%Using formulas (\ref{expected_value}) and (\ref{variance}) 
Using the fact that for independent random variables the expected value of the product is the product of expected values, one can easily show that
\begin{equation}
    \mathbb{E}(Z) = 0,~~ Var(Z) = \frac{n_Y}{n_Y-2} \sigma_X^2,
\end{equation}
when $n_Y>2$.
The expectation and the variance of the product are simply the products of the corresponding moments of the marginal variables. Note that the variance decreases to $\sigma_X^2$ with $n_Y\rightarrow\infty$. 

\begin{figure}[ht]
    \centering
    \includegraphics[scale = 0.5]{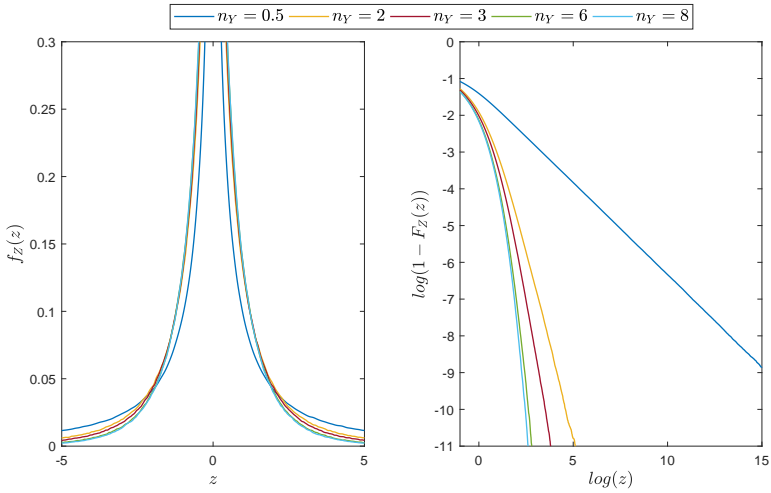}
    \caption{Distribution of the product of two independent random variables from the Gaussian and Student's t distribution with $\mu_X=0$ and $\sigma_X=1$ and different values of $n_Y$ parameter. Left panel: the PDF of the random variable $Z$. Right panel: the distribution tail for the random variable $Z$ (in log-log scale).}
    \label{fig:gauss_student}
\end{figure}

\begin{figure}[ht]
    \centering
    \includegraphics[scale = 0.5]{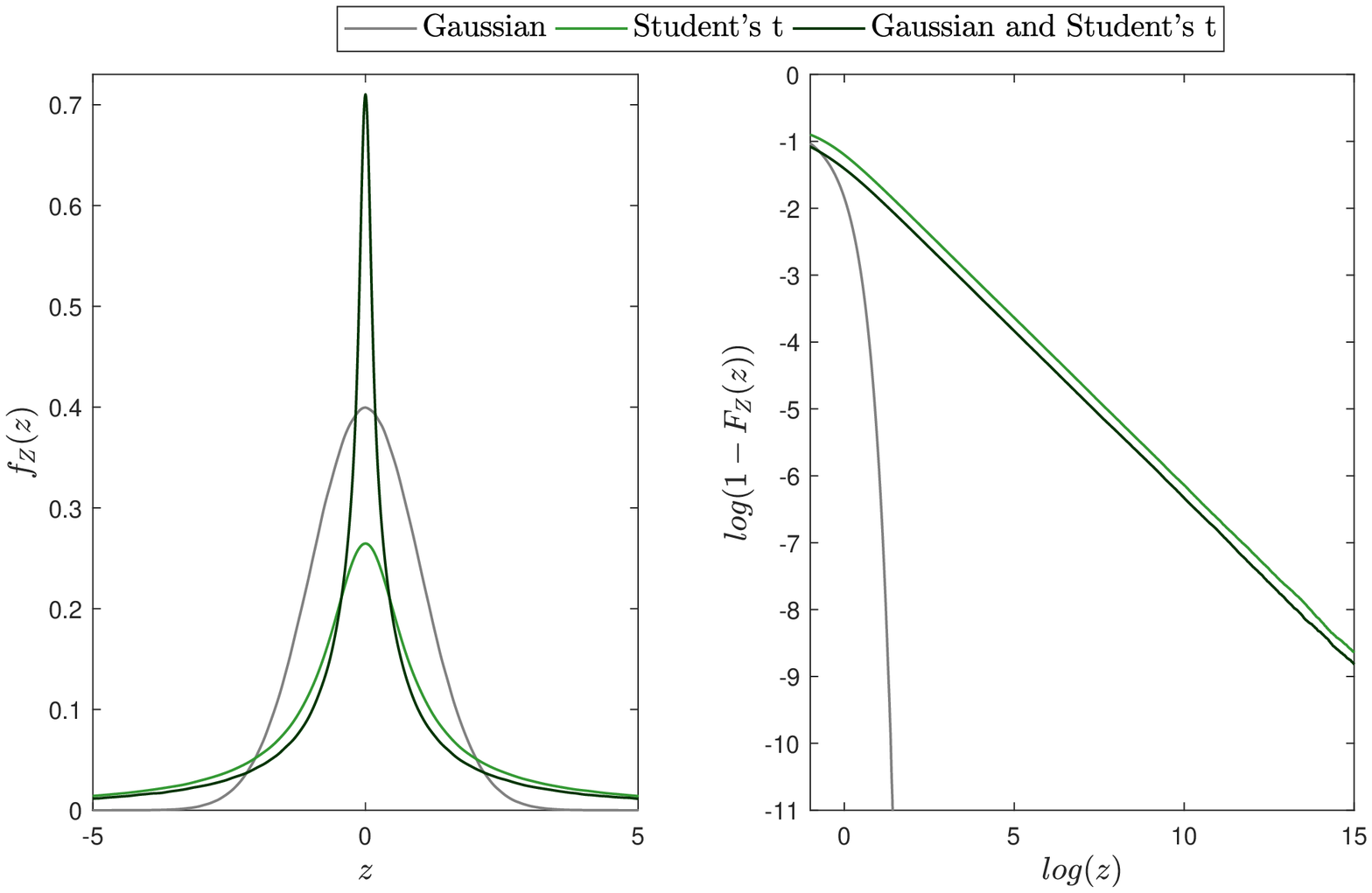}
    \caption{Left panel: the PDF of the random variables $X$ (Gaussian), $Y$ (Student's t) and $Z$. Right panel: the distribution tail's for $X$, $Y$ and $Z$. The parameters are equal to: $\mu_X=0$ and $\sigma_X=1$ and $n_Y=0.5$.}
    \label{fig1_appenddix}
\end{figure}

In Fig.  \ref{fig:gauss_student} we plot the PDF and the corresponding distribution tail of the random variable $Z$ that is a product of  two independent random variables from the standard Gaussian (i.e. $N(0,1)$) and Student's t distributions for different values of the $n_Y$ parameter. %Here we assume that $\mu_X=0$ and $\sigma_X=1$. 
The resulting distribution is symmetric around zero. The shape of its tail is clearly dependent on the value of the $n_Y$ parameter. For larger values of $n_Y$ the tails are close to the Gaussian ones, while for lower $n_Y$ they become much heavier, resembling rather the Student's t tails.
These three distributions, namely, the Gaussian, Student's t and their product, are compared in {Fig. \ref{fig1_appenddix}. One can see that the product has a lighter tail than the Student's t distributed random variable but a heavier tail than the Gaussian one. }

\subsection{Gaussian and Pareto distributions}
%We assume that the one-dimensional random variable $X$ has the Gaussian distribution with parameters $\mu_X\in \mathbb{R}$ and $\sigma_X>0$, while the one-dimensional random variable $Y$ has the Pareto distribution with the PDF given in Eq. (\ref{pareto2}) with $a_Y,\theta_Y>0$. Moreover, for simplicity we consider the case with $\mu_X=0$. %We assume that $X$ and $Y$ are independent, thus \textcolor{blue}{the random vector $(X,Y)$ has the following PDF
%\begin{equation} 
%   f_{X,Y}(x,y) =  \frac{a_Y (\theta_Y)^{a_Y}}{\sqrt{2\pi}\sigma_Xy^{a_Y+1}} \exp\Bigg\{-\frac{x^2}{2\sigma_X^2}\Bigg\},~x\in \mathbb{R}, y>\theta_Y.
%\end{equation}}
%In the following Lemma we give the explicit form of the PDF of the random variable $Z$ being a product of $X$ and $Y$.
Assume that $X$ is a Gaussian distributed random variable with parameters $\mu_X=0$ and $\sigma_X>0$ and $Y$ is a Pareto distributed random variable with parameters $a_Y>0$ and $\theta_Y>0$. Moreover, we assume that $X$ and $Y$ are independent. 
\begin{lemma}\label{lema2} 
The random variable $Z$ defined as a product of $X$ and $Y$ has the following PDF
\begin{eqnarray}
%f_Z(z)=\frac{a_Y \theta_Y^{a_Y}}{2\sqrt{\pi} z^{a_Y+1}} (1-(-1)^{a_Y})
  %(\sqrt{2}\sigma_X)^{a_X}
 % \Gamma\left(\frac{a_Y}{2}+\frac{1}{2}\right),~z\in \mathbb{R},
 f_Z(z) =\frac{a_Y\theta_Y^{a_Y}\sigma_X^{a_Y}}{\sqrt{2\pi}\mid z\mid^{a_Y+1}}2^{(a_Y-1)/2}\gamma\left(\frac{a_Y+1}{2},\frac{z^2}{2\sigma^2\theta_Y^2}\right),~~z\in\mathbb{R},
\end{eqnarray}
where $\Gamma(\cdot)$ is the Gamma function and $\gamma(\cdot,\cdot)$ is the upper incomplete Gamma function defined as $\gamma(\alpha,x)=\int_{0}^{x}t^{\alpha-1}\exp\{-t\}dt$.
\end{lemma}
\noindent The proof of this lemma is given in Appendix A.

Using the fact that $\gamma(a,x)\rightarrow \Gamma(a)$ when $x\rightarrow \infty$, for $z\rightarrow \infty$ we have 
\begin{eqnarray}\label{asympt}
%f_Z(z)=\frac{a_Y \theta_Y^{a_Y}}{2\sqrt{\pi} z^{a_Y+1}} (1-(-1)^{a_Y})
  %(\sqrt{2}\sigma_X)^{a_X}
 % \Gamma\left(\frac{a_Y}{2}+\frac{1}{2}\right),~z\in \mathbb{R},
 f_Z(z)\sim \frac{a_Y\theta_Y^{a_Y}\sigma_X^{a_Y}}{\sqrt{2\pi}\mid z\mid^{a_Y+1}}2^{(a_Y-1)/2}\Gamma\left(\frac{a_y+1}{2}\right). 
 \end{eqnarray}
 This indicates that the PDF of $Z$ has power-law behavior and thus in the considered case the Pareto distribution dominates the tail behavior.

%From formula (\ref{expected_value}) we obtain that 
The expected value of $Z$ in the considered case, when $a_Y>1$, is equal to zero and %from Eq. (\ref{variance}) we get the following formula for 
the variance of $Z$, if $a_Y>2$ is given by
\begin{equation}
      Var(Z) = \frac{ (\theta_Y\sigma_X)^2a_Y}{(a_Y-2)}.
    \end{equation}
The variance decreases to the product of the marginal variables scale parameters squares, $(\theta_Y\sigma_X)^2$, as the Pareto shape parameter $a_Y\rightarrow\infty$.

\begin{figure}[ht]
    \centering
    \includegraphics[scale = 0.5]{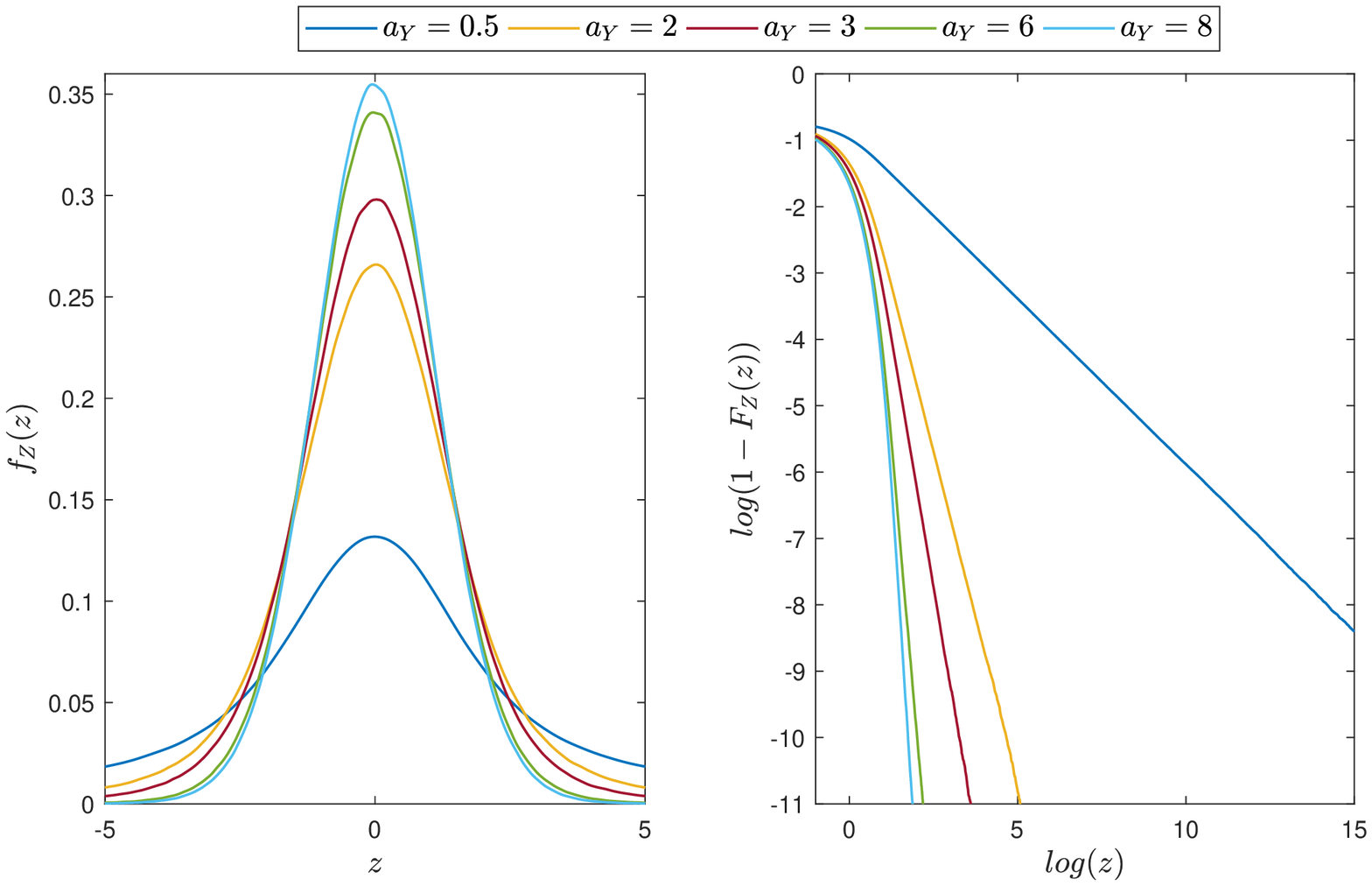}
    \caption{Distribution of the product of two independent random variables from the Gaussian and Pareto distribution with $\mu_X=0$ and $\sigma_X=\theta_Y=1$ and different values of the $a_Y$ parameter. Left panel: the PDF of the random variable $Z$. Right panel: the distribution tail for the random variable $Z$ (in log-log scale).}
    \label{fig:gauss_pareto}
\end{figure}

\begin{figure}[ht]
    \centering
    \includegraphics[scale = 0.5]{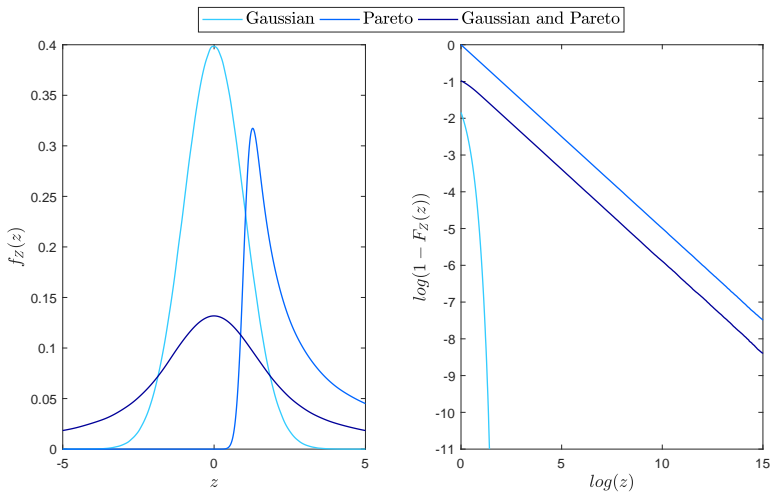}
    \caption{Left panel: the PDF of the random variables $X$ (Gaussian), $Y$ (Pareto) and $Z$. Right panel: the distribution tails. The parameters are: $\mu_X=0$, $\sigma_X=\theta_Y=1$ and $a_Y=0.5$.}
    \label{fig2_appenddix}
\end{figure}

In Fig.  \ref{fig:gauss_pareto}  we plot the PDF and the corresponding distribution tail of the random variable $Z$ that is a product of  two independent random variables from the Gaussian and Pareto distributions for different values of $a_Y$ parameter. Moreover 
%Here we assume that $\mu_X=0$ and $\sigma_X=\theta_Y=1$. 
in Fig. \ref{fig2_appenddix} we demonstrate the comparison of the PDFs and distribution tails of $X$ and $Y$,  and the corresponding random variable $Z$ for selected values of the parameters. One can see that the product has a lighter tail than the Pareto distributed random variable but a heavier tail than the Gaussian one. However, the power-law behavior, which corresponds to  Eq. (\ref{asympt}), can be  easily observed.

\subsection{Log-normal and Student’s t distributions}
In this section, we assume that the one-dimensional random variable $X$ has the log-normal distribution with parameters $\mu_X\in \mathbb{R}$ and $\sigma_X>0$, while the one-dimensional random variable $Y$ has the Student's t distribution with the PDF given in Eq. (\ref{student2}) with $n_Y>0$.  We assume that $X$ and $Y$ are independent.
%, \textcolor{blue}{thus the random vector $(X,Y)$ has the following PDF
%\begin{equation} \label{lognormal_Student}
%   f_{X,Y}(x,y) = \frac{\Gamma{\left(\frac{n_Y+1}{2}\right)}}{\pi x\sigma_X\sqrt{2n_Y}\Gamma\left(\frac{n_Y}{2}\right)} \exp\Bigg\{-\frac{(\log(x)-\mu_X)^2}{2\sigma_X^2}\Bigg\}\left( 1+\frac{y^2}{n_Y} \right)^{-\frac{n_Y+1}{2}},
%\end{equation}
%where $x>0,~y\in \mathbb{R}$.} 
Using Eq.  (\ref{pdf_Z_independent}) one obtains that the PDF of the random variable $Z$ is given by
\begin{align}\label{LogSTU}
f_Z(z)&=\frac{\Gamma{\left(\frac{n_Y+1}{2}\right)}}{\pi \sigma_X\sqrt{2n_Y}\Gamma\left(\frac{n_Y}{2}\right)}\\
&\times \int_{0}^{\infty} \frac{1}{x^2}\exp\Bigg\{-\frac{(\log(x)-\mu_X)^2}{2\sigma_X^2}\Bigg\}\left( 1+\frac{z^2}{n_Yx^2}  \right)^{-\frac{n_Y+1}{2}}dx.\nonumber
\end{align}
The PDF given in Eq. (\ref{LogSTU}) has no closed form representation and requires numerical calculations.
%
%Using formulas (\ref{expected_value}) and (\ref{variance}) 
However, using the independence assumption, one can easily show that
\begin{equation}
    \mathbb{E}(Z) = 0,~~ Var(Z) = \frac{n_Y}{n_Y-2} \exp\Bigg\{\sigma_X^2-1\Bigg\}\exp\Bigg\{2\mu_X+\sigma_X^2\Bigg\},
\end{equation}
when $n_Y>2$.
Similarly as for the Gaussian-Student's t case, the variance decreases to the variance of $X$ with $n_Y\rightarrow\infty$.

\begin{figure}[ht]
    \centering
    \includegraphics[scale = 0.5]{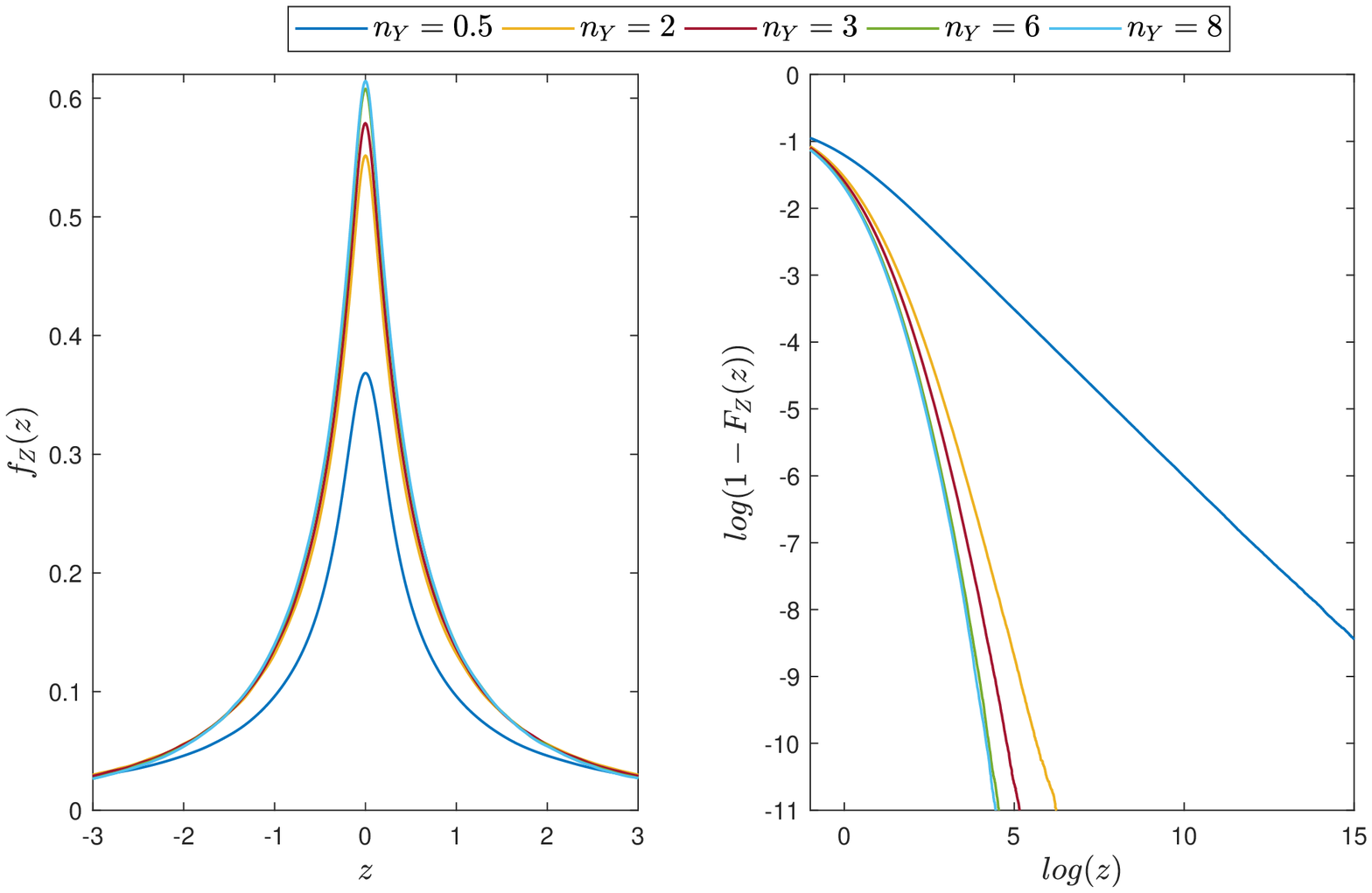}
    \caption{Distribution of the product of two independent random variables from the log-normal and Student's t distribution with $\mu_X=0$ and $\sigma_X=1$ and different values of $n_Y$ parameter. Left panel: the PDF of the random variable $Z$. Right panel: the distribution tail for the random variable $Z$ (in log-log scale). }
    \label{fig:logn_student}
\end{figure}

\begin{figure}[ht]
    \centering
    \includegraphics[scale = 0.5]{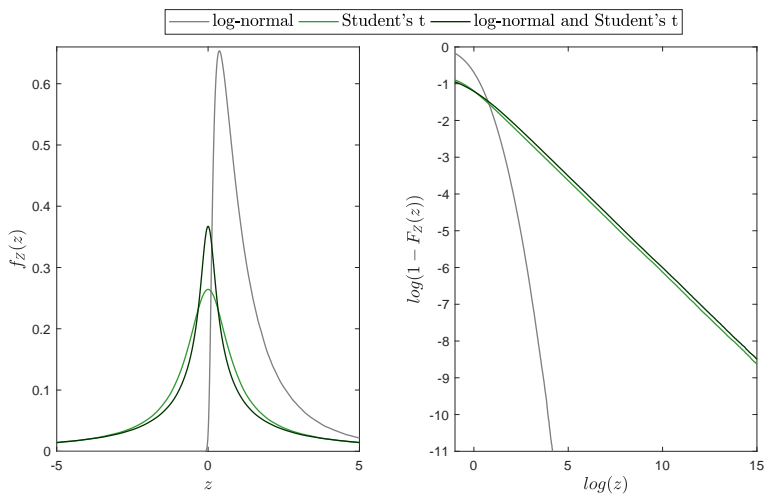}
    \caption{Left panel: the PDF of the random variables $X$ (log-normal), $Y$ (Student's t) and $Z$. Right panel: the distribution tails. The parameters are: $\mu_X=0$ and $\sigma_X=1$ and $n_Y=0.5$.}
    \label{fig3_appenddix}
\end{figure}

In Fig.  \ref{fig:logn_student} we plot the PDF and the corresponding  distribution tail of the random variable $Z$  for different values of the $n_Y$ parameter, while in Fig. \ref{fig3_appenddix} we show a comparison of these three distributions. The picture is similar as the one obtained for the Gaussian - Student's t case (see Fig. \ref{fig:gauss_student}) with slightly heavier tails, especially for higher values of $n_Y$ parameters.  

\subsection{Log-normal and Pareto distributions}
As the last case, %we consider the one-dimensional random variable $X$ having the log-normal distribution with parameters $\mu_X\in \mathbb{R}$ and $\sigma_X>0$, while the one-dimensional random variable $Y$ has the Pareto distribution with the PDF given in Eq. (\ref{pareto2}) with $a_Y,\theta_Y>0$. For simplicity, we consider the case with $\mu_X=0$. We assume that $X$ and $Y$ are independent, \textcolor{blue}{, the random vector $(X,Y)$ has the following PDF
%\begin{equation} 
%   f_{X,Y}(x,y) =  \frac{a_Y (\theta_Y)^{a_Y}}{\sqrt{2\pi}x\sigma_Xy^{a_Y+1}} \exp\Bigg\{-\frac{\log(x)^2}{2\sigma_X^2}\Bigg\},~x>0, y>\theta_Y.
%\end{equation}}
assume that $X$ is a log-normally distributed random variable with parameters $\mu_X=0$ and $\sigma_X>0$ and $Y$ is a Pareto distributed random variable with parameters $a_Y>0$ and $\theta_Y>0$. Moreover, assume that $X$ and $Y$ are independent. 
\begin{lemma}\label{lema3}
The random variable $Z$ defined as a product of $X$ and $Y$ has the following PDF
\begin{eqnarray}\label{asympt1}
 f_Z(z) &=& \frac{\theta_Y}{\sigma_Xz}\phi\left(\frac{\log(z/\theta_Y)}{\sigma_X}\right)\\
 &&+\frac{a_Y\theta_Y^{a_y}\exp\Bigg\{\frac{a_y^2\sigma_X^2}{2}\Bigg\}}{z^{a^y+1}}\Phi\left(\frac{\log(z/\theta_Y)-a_Y\sigma_X^2}{\sigma_X}\right)\nonumber\\
 &&-\frac{\theta_Y^{a_y+1}\exp\Bigg\{\frac{a_y^2\sigma_X^2}{2}\Bigg\}}{\sigma_Xz^{a^y+1}}\phi\left(\frac{\log(z/\theta_Y)-a_Y\sigma_X^2}{\sigma_X}\right),~z>0,\nonumber
\end{eqnarray}
where $\phi(\cdot)$ and $\Phi(\cdot)$ are the PDF and CDF of the standard Gaussian distribution, respectively. 
\end{lemma}
\noindent The proof of this lemma is presented in Appendix B.

The expected value of $Z$ (when $a_Y>1$) and the variance of $Z$ (when $a_Y>2$) are given by
\begin{equation}
\mathbb{E}(Z)=\frac{\theta_Ya_Y}{a_Y-1}\exp\Bigg\{\frac{\sigma_X^2}{2}\Bigg\},~~
      Var(Z) = \frac{\theta_Ya_Y}{(a_Y-2)}\left(\exp\Bigg\{\sigma_X^2\Bigg\}-1\right)\exp\Bigg\{\sigma_X^2\Bigg\}.\nonumber
    \end{equation}
Here, both marginal variables can take only positive values. Hence, the product variable $Z$ is also positive and its expected value (if it exists) is greater than 0, $\mathbb{E}(Z)>0$. The expected value, as well as, the variance of $Z$ are the products of the corresponding moments of the marginal distributions. Their existence is directly related to the existence of the moments of the Pareto marginal distribution.     
    
\begin{figure}[ht]
    \centering
    \includegraphics[scale = 0.5]{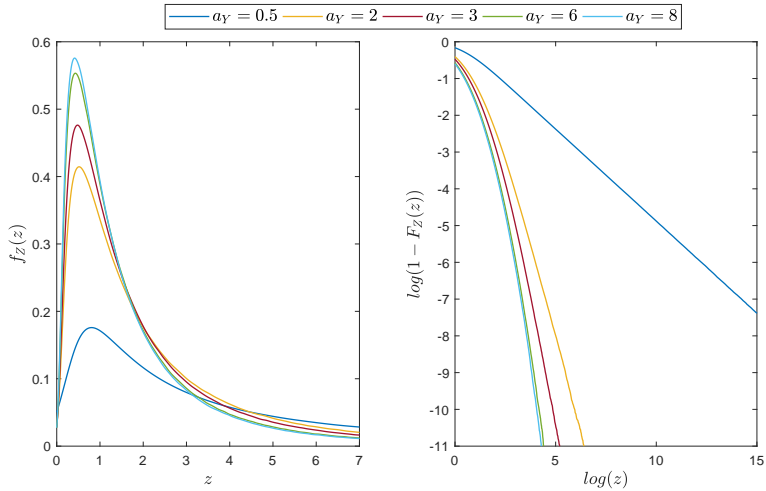}
    \caption{Distribution of the product of two independent random variables from the log-normal and Pareto distribution with $\mu_X=0$ and $\sigma_X=\theta_Y=1$ and different values of $a_Y$ parameter. Left panel: the PDF of the random variable $Z$. Right panel: the distribution tail for the random variable $Z$ (in log-log scale). }
    \label{fig:logn_pareto}
\end{figure}

\begin{figure}[ht]
    \centering
    \includegraphics[scale = 0.5]{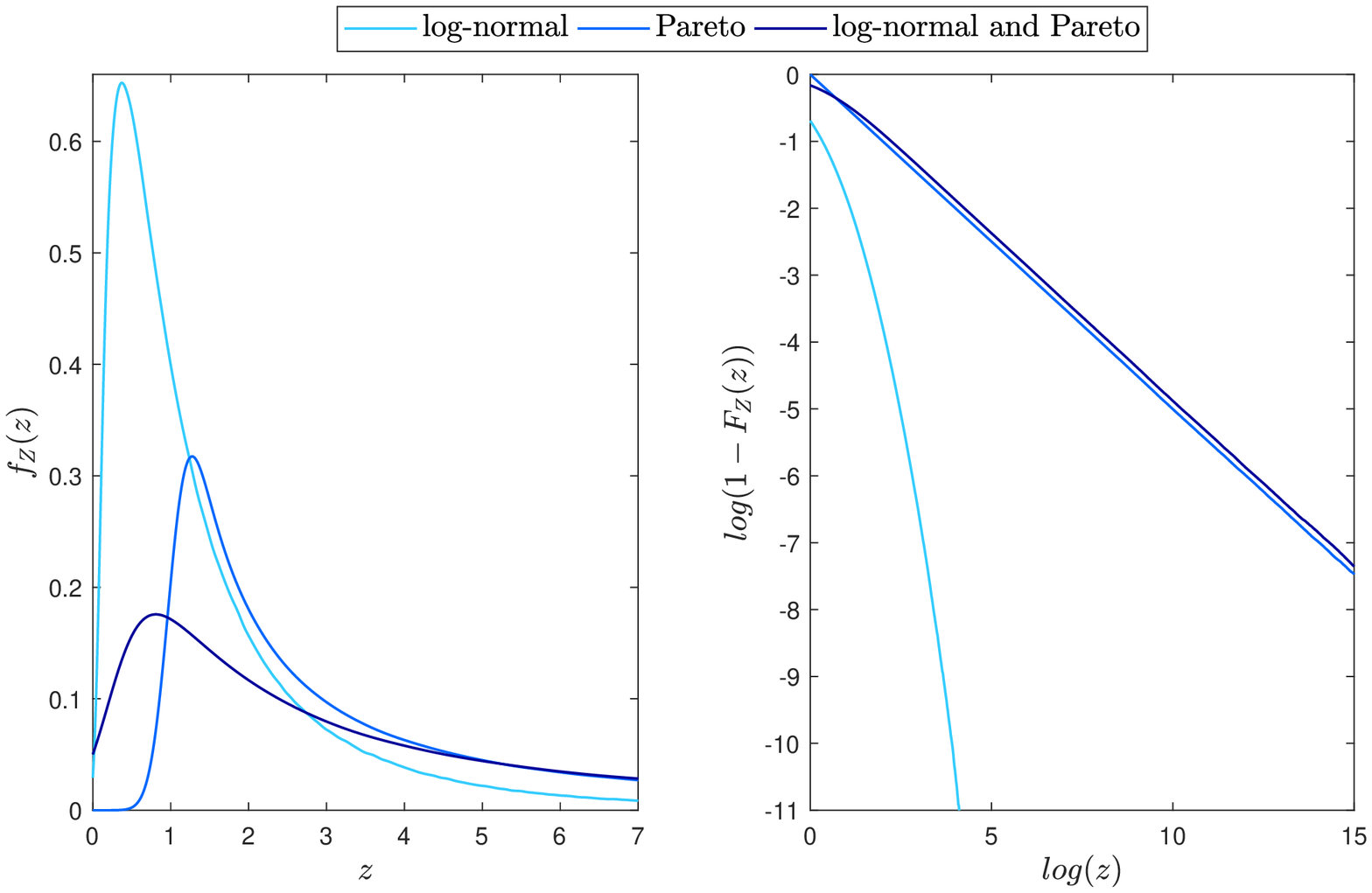}
    \caption{Left panel: the PDF of the random variables $X$ (log-normal), $Y$ (Pareto) and $Z$. Right panel: the distribution tails. The parameters are equal: $\mu_X=0$ and $\sigma_X=\theta_Y=1$ and $a_Y=0.5$.}
    \label{fig4_appenddix}
\end{figure}

In Fig.  \ref{fig:logn_pareto} we plot the PDF of the random variable $Z$ and the corresponding distribution tail for different values of $a_Y$ parameter, 
while in Fig. \ref{fig4_appenddix} we demonstrate the comparison of the $X$, $Y$ and $Z$ distributions. As can be observed the tail parameter of the marginal Pareto distribution has a large impact on the tail of the product distribution. 

\section{Parameters estimation - simulation study}\label{sec:5}
Based on the formulas for the PDF, $f_Z(z)$, derived in the previous sections, we illustrate how the parameters of the  product distribution can be estimated. To this end, we simulate 
samples of random vectors $((X_1,Y_1),(X_2,Y_2),...,(X_N,Y_N))$ from the analyzed distributions and calculate the product of each simulated pair, $Z_i=X_iY_i, i=1,2,...,N$. The Gaussian random vectors are generated using the Cholesky decomposition. The simulation of the vectors from the log-normal and dependent Student's t distribution is based on the relation with the Gaussian distribution. In the case of dependent Pareto variables, the method based on conditional distributions is used, based on the fact that $F_{X,Y}(x,y)=F_X(x)F_{Y\mid X}(y\mid x)$. First, $X$ is generated using the marginal distribution $F_X(x)$. Next, $Y$ is generated using the conditional distribution $F_{Y\mid X}(y\mid x=X)$ with $X$ from the first step. For both steps, the inverse CDF method is applied. The independent Student's t, independent Pareto, Gaussian-Student's t and Gaussian-Pareto vectors simulation is based on the corresponding one-dimensional distributions. 

The estimation of parameters is based on the maximum likelihood method. Thus, we have
\begin{equation}
\hat{\Theta}=\text{argmax}_{\Theta} L(Z_1,Z_2,....,Z_N;\Theta) =\text{argmax}_{\Theta} \prod_{i=1}^N f_Z(Z_i;\Theta), 
\end{equation}
where $\Theta$ is the vector of parameters and $f_Z(z)$ is the corresponding product probability density function. Since analytical solutions do not exist in any of the considered cases, the maximization is performed numerically. Note that for the product of the multivariate Gaussian distribution with $\mu_X=\mu_Y=0$ the method of moments estimators can be also straightforwardly derived, using the exact formulas for the moments.
%Eq. (\ref{Gaussian_moments}). 
Note also that the scale parameters of the individual marginal distributions can not be inferred separately from the product, as they jointly yield the scale of the resulting one-dimensional variable, e.g., for the Gaussian distribution we have $\sigma_X\sigma_Y$, for the log-normal distribution $\sigma_X^2+\sigma_Y^2+2\rho\sigma_X\sigma_Y$, while for the independent Pareto $\theta_X\theta_Y$.

\begin{figure}[ht]
    \centering
    \includegraphics[scale=0.8]{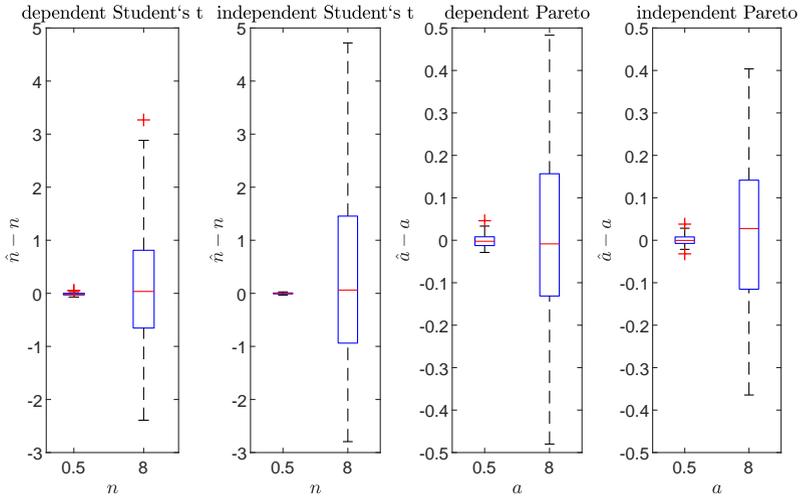}
    \caption{Boxplots of the errors of the degrees of freedom estimator, $\hat{n}$ and the tail parameter, $\hat{a}$, for the Student's t and Pareto distribution products, respectively. Left panels: results for the dependent and independent  Student's t distribution with $\rho=0$ and $n=0.5$ or $n=8$. Right panels: results for the dependent  and independent Pareto  distribution with $\theta_X=\theta_Y=1$ and $a=a_X=a_Y=0.5$ or $a=a_X=a_Y=8$. The length of each sample was equal to 1000, while the number of repetitions was equal to 100.}
    \label{fig:fig2_est}
\end{figure}

In Fig. \ref{fig:fig2_est} we plot the boxplots of the errors of the shape parameter estimation for the product of the  Student's t (i.e., $\hat{n}-n$) and Pareto (i.e. $\hat{a}-a$) distributions. The dependent variables case (see Eq. (\ref{student2}) for Student's t and Eq. (\ref{pareto2}) for Pareto) as well as the independent one (see Eq. (\ref{student4}) for Student's t and Lemma 1 for Pareto) are considered.  The other parameters are set to $\rho=0$, and  $\theta_X=\theta_Y=1$. In the independent Pareto case, the same shape parameter for both variables, $a_X=a_Y=a$, is considered. 
The boxplots are plotted for two simulated values of the shape parameter, namely, $n=0.5$ or $n=8$ and $a=0.5$ or $a=8$. The obtained error distribution is spread around 0, with much larger deviations in the lighter tails case (i.e., $n=8$ and $a=8$). Such an effect might be caused by the fact that for larger values of the shape parameters the differences between distributions become smaller. In the case of the Student's t distribution, the errors are not symmetric around zero with more cases of overestimation of the parameter $n$. In the case of the Pareto distribution, there is no visible asymmetry.

\begin{figure}[ht]
    \centering
    \includegraphics[scale=0.8]{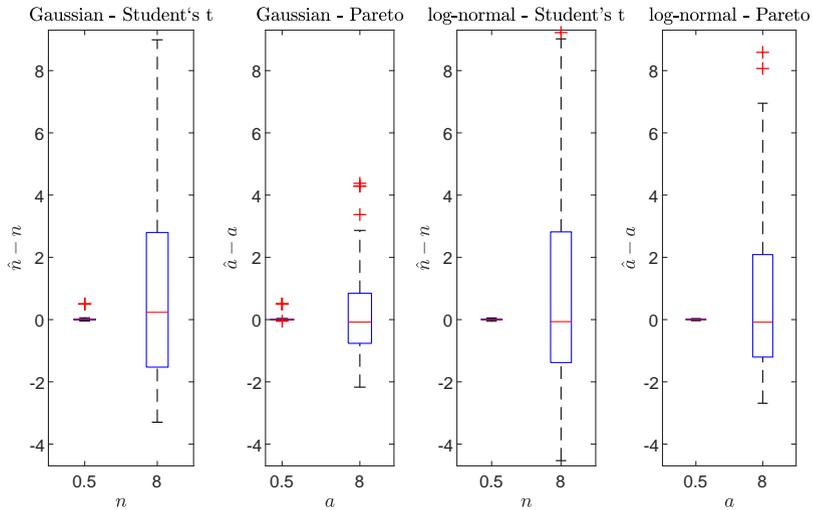}
    \caption{Boxplots of the errors of the degrees of freedom $\hat{n}$ and the tail parameter $\hat{a}$ estimators for the products of the Gaussian or log-normal distribution and Student's t or Pareto distribution, respectively. The parameters of the Gaussian and log-normal distribution were set to $\mu_X=0$ and $\sigma_X=1$. The Student's t distribution was simulated with $n=0.5$ or $n=8$ degrees of freedom, while the Pareto distribution with $\theta=1$ and $a=0.5$ or $a=8$. The length of each sample was equal to 1000, while the number of repetitions was equal to 100.}
    \label{fig:fig3_est}
\end{figure}

Finally, in Fig. \ref{fig:fig3_est} we plot the boxplots of the errors of the shape parameters for the products of different marginal distributions, namely, the Gaussian - Student's t, the Gaussian - Pareto, the log-normal - Student's t and the log-normal - Pareto case. For the Gaussian and log-normal distributions we assume $\mu_X=0$ and $\sigma_X=1$, while the scale parameter of the Pareto distribution is equal $\theta_Y=1$. The errors are analyzed in two cases representing the infinite- ($n=0.5$ for Student's t distribution or $a_Y=0.5$ for Pareto distribution) as well as finite-variance ($n=8$ for Student's t distribution or $a_Y=8$ for Pareto distribution).  The obtained estimation results are similar to the purely Student's t or Pareto distributions (see Fig. \ref{fig:fig2_est}). The higher errors are obtained for the higher values of the $n$ and $a$ parameters. Furthermore, the effect of overestimation of the Student's t degrees of freedom can be noticed.  

\section{Real data application - distribution of electricity transaction values}\label{sec:6}
In this section, we demonstrate the possible application of the theoretical results discussed in the previous sections to a real-data case. We use the transactions data from the German electricity market settled in the EPEX energy exchange. Each transaction is characterized by two values: the volume of sold energy (in MWh) and the price (in EUR/MWh). The final transaction value, being the amount of the total profit for the energy seller or a total cost for energy buyer, is the product of these two variables. Hence, knowing the product distribution is important for profit/costs planning in energy companies. 

The data comes from a continuous trading on the intraday market. A transaction is settled each time two bid and sell offers meet. Hence, each data point corresponds to a different offer and usually a different market participant, so the sample points can be assumed independent. One of the main characteristics of the electricity market is its seasonality on the yearly, weekly, and daily level \cite{RWeron_energy}. To avoid the possible influence of the transaction time on its distribution, we analyze separately the transactions being settled in different hours. The data points within one hour are assumed to be identically distributed. 

We analyze vectors $\{(x_1,y_1),(x_2,y_2),...,(x_n,y_n) \}$, where $x_i$ is the volume and $y_i$ is the price of the $i$-th transaction within a given hour and $n$ is the number of transactions during this hour. We also analyze a sample of the transaction values being a product of each $x_i$ and $y_i$, i.e. $\{x_1 y_1,x_2 y_2,...,x_n y_n \}$. For an illustration, we have chosen two representative hours, namely, hour 14 from 20th December 2020 and hour 8 from 15th March 2020. The first case shows a typical picture that is observable for most hours, while the second case illustrates the less frequent but typical for the German electricity market situation where some of the transaction prices are negative. Negative prices are a unique feature of electricity markets caused by the limited ability to store electricity and the fast development of production from renewable energy sources \cite{NegativePrices}. If there is a short-term oversupply of electricity due, for example, to wind or solar production, the producers from conventional energy sources might be willing to pay for the reception of electricity instead of stopping the production units.  

\begin{figure}[ht]
    \centering
    \includegraphics[scale=0.8]{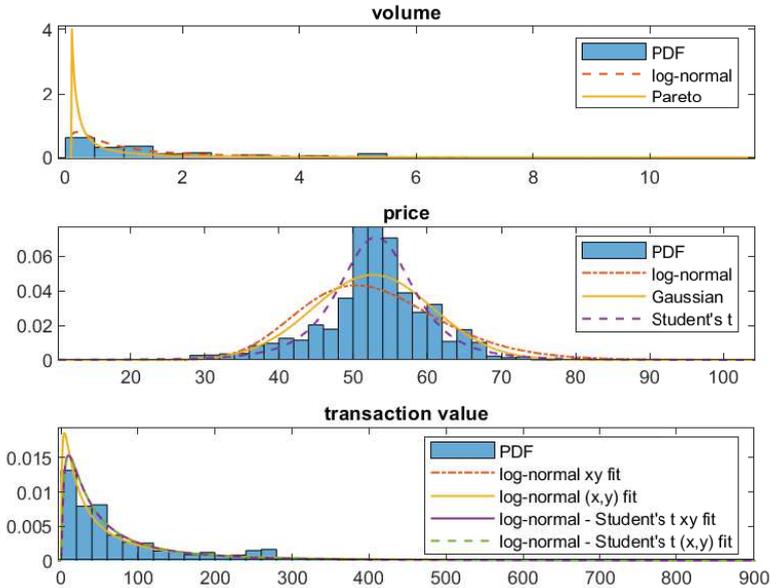}
    \caption{Probability distribution of the volume (top panel), price (middle panel) and their product, i.e. the transaction value (bottom panel) from the continuous intraday energy market data on 20th December 2020, hour 14. The fitted PDFs are also plotted. The transaction value densities are fitted in two ways: separately to the $x$ (price) and $y$ (volume) data (denoted by '$(x,y)$ fit') or to the product $xy$ (transaction value) data (denoted by '$xy$ fit').}
    \label{fig:data1}
\end{figure}

\begin{figure}[ht]
    \centering
    \includegraphics[scale=0.8]{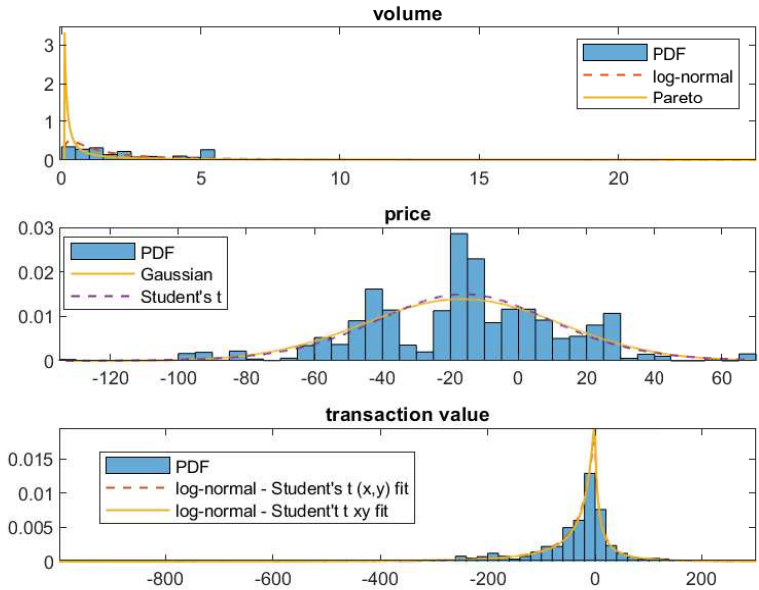}
    \caption{Probability distribution of the volume (top panel), price (middle panel) and their product, i.e. the transaction value (bottom panel) from the continuous intraday energy market data on 15th March 2020, hour 8. The fitted PDFs are also plotted.  The transaction value densities are fitted in two ways: separately to the $x$ (price) and $y$ (volume) data (denoted by '$(x,y)$ fit') or to the product $xy$ (transaction value) data (denoted by '$xy$ fit').}
    \label{fig:data2}
\end{figure}

The corresponding PDFs are plotted in Figs. \ref{fig:data1} and \ref{fig:data2} for 20th December and 15th March, respectively. The correlation between the volume and price samples is not significant. Hence, we assume that the  corresponding $X$ and $Y$ random variables are not correlated. Since the volume can not be negative, the only distributions from the ones analyzed in the paper that can be fitted to the volume sample are the log-normal and Pareto. For the price sample, the Gaussian, log-normal (for positive data) and Student't t with location and scale parameters are fitted. The resulting densities are plotted in the corresponding panels of Figs. \ref{fig:data1} and \ref{fig:data2}. 


For both analyzed hours, the log-normal distribution density function resembles the shape of the volume histogram. On the other hand, the Pareto distribution clearly overestimates the probability of small volumes. Hence, for the product analysis, we will only consider the log-normal distribution for the first coordinate. 
Looking at the price distribution for the data from hour 14 on 20th December 2020 (see Fig. \ref{fig:data1}), we can observe that only the Student's t PDF resembles the shape of the obtained histogram, especially around the mean, where the Gaussian and log-normal probabilities are underestimated. The picture for the second considered data set, i.e., the price distribution from hour 8 on 15th March 2020 (see Fig. \ref{fig:data2}), is different. First, there is a significant probability of obtaining negative values. Hence, in this case we do not fit the log-normal distribution, which can take only positive values.  Second, the Gaussian and the Student's t PDF produce a similar fit. Recall that with increasing degrees of freedom in the latter case, the distribution tends to the Gaussian one. The tails for this hour are lighter than for the 20th December data.  Overall, from the set of considered distributions, only the Student's t can be well fitted to the prices from both analyzed hours.      

Next, we analyze the resulting product, i.e., the transaction value distribution. To this end, we proceed in two ways. Firstly, we use the parameters estimated for the one-dimensional samples of the random variables $X$ and $Y$ (i.e., volumes and prices) and use the formulas for the product PDF derived in Sections \ref{sec:2}-\ref{sec:4}. Based on the results obtained for the one-dimensional samples we show only the product of log-normal distribution or log-normal and Student's t distribution (log-normal-Student's t).
Secondly, we also fit the derived PDFs to the sample of the product $XY$.
 We can observe that the final product PDF obtained by these two ways in the case of log-normal-Student's t overlaps and resembles the shape of the sample histogram, confirming a good fit. As mentioned before,  due to the negative prices apparent in the second analyzed hour, the log-normal distribution can be fitted only to the first dataset. In this case, we can observe a discrepancy between the density obtained using the $(X,Y)$ coordinates and its product $XY$.  Comparing the obtained shapes of the distribution for both hours, we can see a clear difference between these cases. The transaction values for hour 8 on 15th March 2020 take both negative and positive values, with long left tails of the distribution. On the other hand, there are no negative transaction values for hour 14 on 20th December 2020 and the distribution is right skewed. Overall, only the log-normal-Student's t case from the considered distributions set is flexible enough to be reasonably applied to both analyzed representative hours and the obtained distribution resembles the data well.

\section{Conclusions}
In this paper, we have discussed the distribution of a random variable that is a product of two continuous random variables. The main attention was paid to the influence of the parameters of the marginal random variables on the final product characteristics. We have considered exemplary distributions belonging to the classes of finite as well as infinite-variance distributions. For both cases, we have discussed how the correlation coefficient between the marginal random variables influences the probabilistic properties of the product. In the case of the Gaussian and log-normal distributions, the non-zero correlation coefficient indicates the statistical dependence of the random variables, while in the second case, i.e., the Student's t and Pareto distributions, this statement is not true.  The most interesting case was the product of the random variables coming from different classes of distributions, for which we have derived the explicit form of the PDF for the Gaussian-Pareto and the log-normal-Pareto case. 

The theoretical results were applied in the proposed estimation methodology. We used a general approach based on the maximum likelihood technique. The presented Monte Carlo simulations clearly indicate the effectiveness of the algorithm. Finally, the  real data analysis was presented. Based on the data from continuous trading on the German energy market, we have shown a good reasonable fit of the product of log-normal and Student's t distribution to the transaction values. Since the transaction value is the final profit/cost for a trader, finding a proper density describing its distribution, which is also consistent with the prices and volumes data, can help an energy market participant in strategy planning.   

Another potential application of the proposed methodology could be the market risk area in metals and mining business.  Mining companies are exposed to two or more market risk factors, like metal prices and currency exchange rates.  These factors behavior is often characterized by non-Gaussian distributions, 
what can be reflected by the discussed in this paper, log-normal and Student’s t distributions. 
From a business perspective of an international mining company, which trades the excavated resources in different than national currency, it is valuable to analyse the selling price also in the national currency, that is a product of the commodity price (usually quoted in USD) and USD/national currency exchange rate. A case study related to the modelling of the copper prices for a polish mining company one can find in \cite{reso}. It has been shown that the behavior of metal price in PLN exhibits specific characteristics, that follow directly from the properties of the individual variables. These properties need to be reflected in optimizing strategies aimed at mitigation of the unacceptable for company market risk. 

\backmatter

%\bmhead{Supplementary information}

%If your article has accompanying supplementary file/s please state so here. 

%Authors reporting data from electrophoretic gels and blots should supply the full unprocessed scans for key as part of their Supplementary information. This may be requested by the editorial team/s if it is missing.

%Please refer to Journal-level guidance for any specific requirements.

\section*{Acknowledgments}
The work of A.W. was supported by National Center of Science under Opus Grant 2020/37/B/HS4/00120 "Market risk model identification and validation using novel statistical, probabilistic, and machine learning tools". J.J. acknowledges a support of NCN Sonata Grant No. 2019/35/D/HS4/00369.

\section*{Data Availability}
The datasets generated during the current study are available from the corresponding author on reasonable request. The energy transactions data analysed in this study are available from 
%the corresponding author upon reasonable request and with permission of 
the EPEX SPOT exchange.

\bibliography{mybibliographysp}

%\section*{Declarations}

%Some journals require declarations to be submitted in a standardised format. Please check the Instructions for Authors of the journal to which you are submitting to see if you need to complete this section. If yes, your manuscript must contain the following sections under the heading `Declarations':

%\begin{itemize}
%\item Funding: The work of A.W. was supported by National Center of Science under Opus Grant 2020/37/B/HS4/00120 "Market risk model identification and validation using novel statistical, probabilistic, and machine learning tools". J.J. acknowledges a support of NCN Sonata Grant No. 2019/35/D/HS4/00369. 

%\item Conflict of interest/Competing interests: The authors have no conflicts of interest to declare that are relevant to the content of this article.
%\item Ethics approval: 
%\item Consent for publication
%\item Availability of data and materials: The datasets generated during and/or analysed during the current study are available from the corresponding author on reasonable request.
%\item Authors' contributions
%\end{itemize}

\noindent
%If any of the sections are not relevant to your manuscript, please include the heading and write `Not applicable' for that section. 

\begin{appendices}

\section{}
\noindent
\textbf{Proof of Lemma \ref{lema2}.} \\
%Using formula (\ref{pdf_Z_independent}) we obtain
%\begin{eqnarray*} 
 %  f_Z(z) &=& \frac{a_Y \theta_Y^{a_Y}}{\sqrt{2\pi}\sigma_X} \int_{-\infty}^{\infty} \frac{1}{|x|} \left(\frac{x}{z}\right)^{a_Y+1} \exp\Bigg\{-\frac{x^2}{2\sigma_X^2}\Bigg\} dx \\
%&=&\frac{a_Y \theta_Y^{a_Y}}{\sqrt{2\pi}\sigma_X z^{a_Y+1}}
 %  \left(\int_{0}^{\infty} x^{a_Y} \exp\Bigg\{-\frac{x^2}{2\sigma_X^2}\Bigg\} dx
  % - \int_{-\infty}^{0} x^{a_Y} \exp\Bigg\{-\frac{x^2}{2\sigma_X^2}\Bigg\} dx \right)\\
   %&=&\frac{a_Y (\theta_Y)^{a_Y}}{\sqrt{2\pi}\sigma_X z^{a_Y+1}} (1-(-1)^{a_Y})
   %\int_{0}^{\infty} x^{a_Y} \exp\Bigg\{-\frac{x^2}{2\sigma_X^2}\Bigg\} dx.
%\end{eqnarray*}
%Taking the substitution $u=\frac{x}{\sqrt{2}\sigma_X}$ in the above integral we get the following
%\begin{eqnarray*}
%\frac{a_Y \theta_Y^{a_Y}}{\sqrt{2\pi}\sigma_X z^{a_Y+1}} (1-(-1)^{a_Y})
 % (\sqrt{2}\sigma_X)^{a_X+1} \int_{0}^{\infty} u^{a_Y} \exp\{-u^2\} du 
  % =\frac{a_Y \theta_Y^{a_Y}}{2\sqrt{\pi} z^{a_Y+1}} (1-(-1)^{a_Y})
  %(\sqrt{2}\sigma_X)^{a_X}
  %\Gamma\left(\frac{a_Y}{2}+\frac{1}{2}\right),
%\end{eqnarray*}
Let us first note that the CDF of $Y$ is given by
\begin{eqnarray}\label{pareto_cdf}
F_Y(y)=1-\theta_Y^{a_Y}y^{-a_Y},~y>\theta_Y.
\end{eqnarray}
Thus, the CDF of $Z$ takes the form
\begin{eqnarray*}
F_Z(z)&=&P(XY<z)=\int_{-\infty}^{\infty}P(Yx<z)f_X(x)dx\\
&=&\frac{1}{\sqrt{2\pi}\sigma_X}\int_{-\infty}^{\infty}P(Yx<z)\exp\Bigg\{-\frac{x^2}{2\sigma_X^2}\Bigg\}dx\\
&=&\frac{1}{\sqrt{2\pi}\sigma_X}\int_{-\infty}^{0}P(Yx<z)\exp\Bigg\{-\frac{x^2}{2\sigma_X^2}\Bigg\}dx\\
&&+\frac{1}{\sqrt{2\pi}\sigma_X}\int_{0}^{\infty}P(Yx<z)\exp\Bigg\{-\frac{x^2}{2\sigma_X^2}\Bigg\}dx\\
&=&\frac{1}{\sqrt{2\pi}\sigma_X}\int_{-\infty}^{0}P\left(Y>\frac{z}{x}\right)\exp\Bigg\{-\frac{x^2}{2\sigma_X^2}\Bigg\}dx\\
&&+\frac{1}{\sqrt{2\pi}\sigma_X}\int_{0}^{\infty}P\left(Y<\frac{z}{x}\right)\exp\Bigg\{-\frac{x^2}{2\sigma_X^2}\Bigg\}dx.
\end{eqnarray*}
We consider separately $z>0$, $z<0$ and $z=0$. For $z>0$ one obtains
\begin{eqnarray*}
F_Z(z)&=&\frac{1}{\sqrt{2\pi}\sigma_X}\int_{-\infty}^{0}\exp\Bigg\{-\frac{x^2}{2\sigma_X^2}\Bigg\}dx\\
&&+\frac{1}{\sqrt{2\pi}\sigma_X}\int_{0}^{\infty}F_Y\left(\frac{z}{x}\right)\exp\Bigg\{-\frac{x^2}{2\sigma_X^2}\Bigg\}dx\\
&=&\frac{1}{\sqrt{2\pi}\sigma_X}\int_{-\infty}^{0}\exp\Bigg\{-\frac{x^2}{2\sigma_X^2}\Bigg\}dx\\&&+\frac{1}{\sqrt{2\pi}\sigma_X}\int_{0}^{z/\theta_Y}\left(1-\theta_Y^{a_Y}\left(\frac{x}{z}\right)^{a_Y}\right)\exp\Bigg\{-\frac{x^2}{2\sigma_X^2}\Bigg\}dx\\
&=&F_X\left(\frac{z}{\theta_Y}\right)-\frac{\theta_Y^{a_Y}}{\sqrt{2\pi}\sigma_Xz^{a_Y}}\int_{0}^{z/\theta_Y}x^{a_Y}\exp\Bigg\{-\frac{x^2}{2\sigma_X^2}\Bigg\}dx.
\end{eqnarray*}
Now, calculating the derivative of $F_Z(z)$ with respect to $z$ for $z>0$ one has 
\begin{eqnarray*}
f_Z(z)&=&\frac{1}{\theta_Y}f_X\left(\frac{z}{\theta_Y}\right)+\frac{a_Y\theta_Y^{a_Y}}{\sqrt{2\pi}\sigma_Xz^{a_Y+1}}\int_{0}^{z/\theta_Y}x^{a_Y}\exp\Bigg\{-\frac{x^2}{2\sigma_X^2}\Bigg\}dx\\
&&-\frac{\theta_Y^{a_Y}}{\sqrt{2\pi}\sigma_Xz^{a_Y}\theta_Y}\left(\frac{z}{\theta_Y}\right)^{a_Y}\exp\Bigg\{-\frac{z^2}{2\sigma_X^2\theta_Y^2}\Bigg\}\\
&=&\frac{1}{\sqrt{2\pi}\sigma_X\theta_Y}\exp\Bigg\{-\frac{z^2}{2\sigma_X^2\theta_Y^2}\Bigg\}\\
&&+\frac{a_Y\theta_Y^{a_Y}}{\sqrt{2\pi}\sigma_Xz^{a_Y+1}}\int_{0}^{z/\theta_Y}x^{a_Y}\exp\Bigg\{-\frac{x^2}{2\sigma_X^2}\Bigg\}dx\\
&&-\frac{1}{\sqrt{2\pi}\sigma_X\theta_Y}\exp\Bigg\{-\frac{z^2}{2\sigma_X^2\theta_Y^2}\Bigg\}\\
&=&\frac{a_Y\theta_Y^{a_Y}}{\sqrt{2\pi}\sigma_Xz^{a_Y+1}}\int_{0}^{z/\theta_Y}x^{a_Y}\exp\Bigg\{-\frac{x^2}{2\sigma_X^2}\Bigg\}dx\\
&=&\frac{a_Y\theta_Y^{a_Y}}{\sqrt{2\pi}\sigma_Xz^{a_Y+1}}2^{(a_Y-1)/2}\\
&& \times\left(\frac{z}{\theta_Y}\right)^{a_Y+1}\left(\frac{\theta_Y\sigma_X}{z}\right)^{a_Y+1}\left(\Gamma\left(\frac{a_Y+1}{2}\right)-\Gamma\left(\frac{a_Y+1}{2},\frac{z^2}{2\sigma^2\theta_Y^2}\right)\right)\\
&=&\frac{a_Y\theta_Y^{a_Y}\sigma_X^{a_Y}}{\sqrt{2\pi}z^{a_Y+1}}2^{(a_Y-1)/2}
\left(\Gamma\left(\frac{a_Y+1}{2}\right)-\Gamma\left(\frac{a_Y+1}{2},\frac{z^2}{2\sigma^2\theta_Y^2}\right)\right)\\
&=&\frac{a_Y\theta_Y^{a_Y}\sigma_X^{a_Y}}{\sqrt{2\pi}z^{a_Y+1}}2^{(a_Y-1)/2}\gamma\left(\frac{a_Y+1}{2},\frac{z^2}{2\sigma^2\theta_Y^2}\right),
\end{eqnarray*}
where $\gamma(a,x)$ is the upper incomplete Gamma function and $\Gamma(a,x)=\Gamma(a)-\gamma(a,x)$ is the lower incomplete Gamma function. For $z<0$ the CDF of $Z$ is given by
\begin{eqnarray*}
F_Z(z)&=&\frac{1}{\sqrt{2\pi}\sigma_X}\int_{-\infty}^{0}P\left(Y>\frac{z}{x}\right)\exp\Bigg\{-\frac{x^2}{2\sigma_X^2}\Bigg\}dx\\
&=&\frac{1}{\sqrt{2\pi}\sigma_X}\int_{-\infty}^{z/\theta_Y}P\left(Y>\frac{z}{x}\right)\exp\Bigg\{-\frac{x^2}{2\sigma_X^2}\Bigg\}dx\\
&&+\frac{1}{\sqrt{2\pi}\sigma_X}\int_{z/\theta_Y}^{0}P\left(Y>\frac{z}{x}\right)\exp\Bigg\{-\frac{x^2}{2\sigma_X^2}\Bigg\}dx\\
&=&F_X\left(\frac{z}{\theta_Y}\right)+\frac{1}{\sqrt{2\pi}\sigma_X}\int_{z/\theta_Y}^{0}\theta_Y^{a_Y}\left(\frac{x}{z}\right)^{a_Y}\exp\Bigg\{-\frac{x^2}{2\sigma_X^2}\Bigg\}dx\\
&=&F_X\left(\frac{z}{\theta_Y}\right)+\frac{\theta_Y^{a_Y}}{\sqrt{2\pi}\sigma_X(-z)^{a_y}}\int_{z/\theta_Y}^{0}(-x)^{a_Y}\exp\Bigg\{-\frac{x^2}{2\sigma_X^2}\Bigg\}dx\\
&=&F_X\left(\frac{z}{\theta_Y}\right)+\frac{\theta_Y^{a_Y}}{\sqrt{2\pi}\sigma_X(-z)^{a_y}}\int^{-z/\theta_Y}_{0}x^{a_Y}\exp\Bigg\{-\frac{x^2}{2\sigma_X^2}\Bigg\}dx.
\end{eqnarray*}
Now, taking the derivative of $F_Z(z)$ with respect to $z$ for $z<0$ we obtain the following formula for corresponding PDF 
\begin{eqnarray*}
f_Z(z)&=&f_X\left(\frac{z}{\theta_Y}\right)-\frac{a_Y\theta_Y^{a_Y}}{\sqrt{2\pi}\sigma_X(-z)^{a_y+1}}\int^{-z/\theta_Y}_{0}x^{a_Y}\exp\Bigg\{-\frac{x^2}{2\sigma_X^2}\Bigg\}dx\\
&&-\frac{\theta_Y^{a_Y}}{\sqrt{2\pi}\sigma_X(-z)^{a_y}\theta_Y}\left(-\frac{z}{\theta_Y}\right)^{a_Y}\exp\Bigg\{-\frac{z^2}{2\sigma_X^2\theta_Y^2}\Bigg\}\\
&=&\frac{1}{\sqrt{2\pi}\sigma_X\theta_Y}\exp\Bigg\{-\frac{z^2}{2\sigma_X^2\theta_Y^2}\Bigg\}\\
&&+\frac{a_Y\theta_Y^{a_Y}}{\sqrt{2\pi}\sigma_X(-z)^{a_y+1}}2^{(a_Y-1)/2}\left(-\frac{z}{\theta_Y}\right)^{a_Y+1}\left(-\frac{\theta_Y\sigma_X}{z}\right)^{a_Y+1}\\
&& \times\left(\Gamma\left(\frac{a_Y+1}{2}\right)-\Gamma\left(\frac{a_Y+1}{2},\frac{z^2}{2\sigma_X^2\theta_Y^2}\right)\right)\\
&&-\frac{1}{\sqrt{2\pi}\sigma_X\theta_Y}\exp\Bigg\{-\frac{z^2}{2\sigma_X^2\theta_Y^2}\Bigg\}\\&=&\frac{a_Y\theta_Y^{a_Y}\sigma_X^{a_Y}}{\sqrt{2\pi}(-z)^{a_Y+1}}2^{(a_Y-1)/2}\left(\Gamma\left(\frac{a_Y+1}{2}\right)-\Gamma\left(\frac{a_Y+1}{2},\frac{z^2}{2\sigma_X^2\theta_Y^2}\right)\right)\\&=&\frac{a_Y\theta_Y^{a_Y}\sigma_X^{a_Y}}{\sqrt{2\pi}(-z)^{a_Y+1}}2^{(a_Y-1)/2}\gamma\left(\frac{a_Y+1}{2},\frac{z^2}{2\sigma_X^2\theta_Y^2}\right).
\end{eqnarray*}
The last case is $z=0$. In this case, the CDF of $Z$ is given by
\begin{eqnarray*}
F_Z(0)=P(XY<0)=P(X<0)=\frac{1}{2}.
\end{eqnarray*}
One can also show that 
\begin{eqnarray*}
\lim_{z\rightarrow 0+}F_Z(z)=\lim_{z\rightarrow 0-}F_Z(z)=\frac{1}{2}=F_Z(0).
\end{eqnarray*}
Thus, finally we obtain the thesis. $\Box$

\section{}
\noindent
\textbf{Proof of Lemma \ref{lema3}.}\\ 
Using the same reasoning as in the proof of Lemma \ref{lema2} one obtains that the CDF of $Z$ takes the following form for $z>0$
\begin{eqnarray*}
F_Z(z)&=&P(XY<z)=\int_{-\infty}^{\infty}P(Yx<z)f_X(x)dx\\
&=&\frac{1}{\sqrt{2\pi}\sigma_X}\int_{0}^{\infty}\frac{1}{x}P(Yx<z)\exp\Bigg\{-\frac{\log(x)^2}{2\sigma_X^2}\Bigg\}dx\\
&=&\frac{1}{\sqrt{2\pi}\sigma_X}\int_{0}^{\infty}\frac{1}{x}P\left(Y<\frac{z}{x}\right)\exp\Bigg\{-\frac{\log(x)^2}{2\sigma_X^2}\Bigg\}dx.
\end{eqnarray*}
Thus, from (\ref{pareto_cdf}) we obtain for $z>0$
\begin{eqnarray*}
F_Z(z)&=&\frac{1}{\sqrt{2\pi}\sigma_X}\int_{0}^{\infty}\frac{1}{x}F_Y\left(\frac{z}{x}\right)\exp\Bigg\{-\frac{\log(x)^2}{2\sigma_X^2}\Bigg\}dx\\
&=&\frac{1}{\sqrt{2\pi}\sigma_X}\int_{0}^{z/\theta_Y}\frac{1}{x}\left(1-\theta_Y^{a_Y}\left(\frac{x}{z}\right)^{a_Y}\right)\exp\Bigg\{-\frac{\log(x)^2}{2\sigma_X^2}\Bigg\}dx\\
&=&\frac{1}{\sqrt{2\pi}\sigma_X}\int_{0}^{z/\theta_Y}\frac{1}{x}\exp\Bigg\{-\frac{\log(x)^2}{2\sigma_X^2}\Bigg\}dx\\
&&-\frac{\theta_Y^{a_y}}{\sqrt{2\pi}\sigma_Xz^{a^y}}\int_{0}^{z/\theta_Y}x^{a_Y-1}\exp\Bigg\{-\frac{\log(x)^2}{2\sigma_X^2}\Bigg\}dx.
\end{eqnarray*}
Now, taking the substitution $\log(x)=u$ in the both above integrals, we obtain
\begin{eqnarray*}
F_Z(z)&=&\frac{1}{\sqrt{2\pi}\sigma_X}\int_{-\infty}^{\log(z/\theta_Y)}\exp\Bigg\{-\frac{u^2}{2\sigma_X^2}\Bigg\}du\\
&&-\frac{\theta_Y^{a_y}}{\sqrt{2\pi}\sigma_Xz^{a^y}}\int_{-\infty}^{\log(z/\theta_Y)}\exp\Bigg\{ua_y-\frac{u^2}{2\sigma_X^2}\Bigg\}du.
\end{eqnarray*}
To make the calculations simpler, in the first integral we put the substitution $w=u/\sigma_X$ and in the second one $w=(u-\sigma_X^2a_Y)/\sigma_X$. Then we obtain
\begin{eqnarray*}
F_Z(z)&=&\int_{-\infty}^{\frac{\log(z/\theta_Y)}{\sigma_X}}\frac{1}{\sqrt{2\pi}}\exp\Bigg\{-\frac{w^2}{2}\Bigg\}dw\\
&&-\frac{\theta_Y^{a_y}\exp\Bigg\{\frac{a_y^2\sigma_X^2}{2}\Bigg\}}{z^{a^y}}\int_{-\infty}^{\frac{\log(z/\theta_Y)-a_Y\sigma_X^2}{\sigma_X}}\frac{1}{\sqrt{2\pi}}\exp\Bigg\{-\frac{w^2}{2}\Bigg\}dw.
\end{eqnarray*}
Calculating the derivative of $F_Z(z)$ with respect to $z>0$ one has 
\begin{eqnarray*}
f_Z(z)&=&\frac{\theta_Y}{\sigma_Xz}\phi\left(\frac{\log(z/\theta_Y)}{\sigma_X}\right)+\frac{a_Y\theta_Y^{a_y}\exp\Bigg\{\frac{a_y^2\sigma_X^2}{2}\Bigg\}}{z^{a^y+1}}\Phi\left(\frac{\log(z/\theta_Y)-a_Y\sigma_X^2}{\sigma_X}\right)\\
&&-\frac{\theta_Y^{a_y+1}\exp\Bigg\{\frac{a_y^2\sigma_X^2}{2}\Bigg\}}{\sigma_Xz^{a^y+1}}\phi\left(\frac{\log(z/\theta_Y)-a_Y\sigma_X^2}{\sigma_X}\right),
\end{eqnarray*}
where $\phi(\cdot)$ and $\Phi(\cdot)$ are the PDF and CDF of the standard Gaussian distribution, respectively.  $\Box$

\end{appendices}



\end{document}